\documentclass[12pt,a4paper]{article}

\usepackage{amsfonts}
\usepackage{amsmath}
\usepackage{amsbsy}
\usepackage{theorem}
\usepackage{graphicx}

\newtheorem{theorem}{Theorem}

\newtheorem{remark}{Remark}

\begin{document}
\sf
\title{\bfseries \large High-level moving excursions for spatiotemporal Gaussian random fields with long range dependence}

\date{}

 \maketitle

\bigskip

\centerline{N.N. Leonenko and M.D. Ruiz--Medina}

\begin{abstract}
The asymptotic behavior of an extended family of integral  geometric random functionals, including
spatiotemporal Minkowski functionals under moving levels, is analyzed in this paper. Specifically, sojourn
measures of  spatiotemporal long--range dependence (LRD)  Gaussian  random fields  are considered in this analysis.
  The limit results derived  provide general reduction principles  under increasing domain asymptotics in space and time. The case of time--varying thresholds is also studied.  Thus, the  family of  morphological measures considered  allows the   statistical and geometrical analysis of random physical systems displaying structural changes over time. Motivated by cosmological applications, the derived results are applied to the context of sojourn
measures of spatiotemporal spherical Gaussian random fields.
 The results are illustrated for  some families of  spatiotemporal  Gaussian  random fields displaying complex spatiotemporal dependence structures.
\end{abstract}


\medskip


\section{Introduction}

\subsection{Connections with Statistical Physics}
\label{sec:1}
Geometric integral  functionals defining  sojourn
measures in the context of random fields play a crucial role in Statistical Physics.
These morphological descriptors arising in integral geometry allow the characterization of connectivity, content and shape of stochastic spatial structures in the analysis of physical systems.  In particular, they  complement  the  spatial statistical and physical analysis of systems, whose stochastic structures are characterized by the Boolean model. Several physical phenomena involved in different  issues of Statistical Physics, such as  complex ﬂuids, porous media, and pattern formation in dissipative systems, can be suitably analyzed in terms of these morphological measures. Particularly, Minkowski functionals, as additive functionals of spatial patterns, provide the geometric description of physical phenomena   from  integral instead of differential expressions.  Special attention has been paid to Minkowski formalism in the three dimensional Euclidean space, where this family of functionals consists of the volume, the surface area of the pattern, its integral mean curvature, and the Euler characteristic or  integral of the Gaussian curvature.

The additivity property of Minkowski functionals is well--suited and applicable in several areas of statistical mechanics where systems have structures or properties that can be decomposed and analyzed in a linear manner.  Some scenarios where additivity of Minkowski functionals is suitable, and hence,  effectively applied are homogeneous and isotropic systems.  Additivity is used in porous media or composite materials to quantify geometric properties like volume, surface area, and connectivity. For well-separated pores or grains, the functional is evaluated over the entire structure by summing  of the functional values  for each component.
This approach is useful in estimating  properties such as permeability, mechanical strength, and conductivity by decomposing the overall geometry into additive components.  In percolation theory, additivity is applicable in the analysis of disconnected clusters.
For systems below the percolation threshold, where clusters are typically separated and finite, the additivity property allows easy computation of geometric measures by summing the contributions of individual clusters.
In materials science, Minkowski functionals are used to describe the geometric features of grains. If the grains are mostly non-overlapping or exhibit simple intersections, the additivity property is useful
 in characterizing grain growth, phase boundaries, and interfaces in such materials, simplifying the calculation of volume, surface area, and curvature. Summarizing, in systems where interactions are weak or the geometry approximates simple shapes (like spheres or cubes), Minkowski functionals can be applied.
Assuming limited or negligible overlap,  in density fields or point distributions, like Poisson point processes or random tessellations, the Minkowski functionals for each region can be summed to describe the global geometry.

Minkowski formalism  offers robust morphological measures contributing to the introduction of order parameters characterizing pattern  transitions in dissipative systems, dynamical quantities characterizing spinodal decomposition, generalized molecular distribution functions for the characterization of the atomic structure of simple fluids (see, eg., Section 3 in \cite{Mecke} and the  references therein). Related statistical physical problems arising in these fields can be solved applying mathematical properties of Minkowski functionals (e.g., computation of virial coefficients and definition of accurate density functionals, prediction of percolation thresholds, and formulation of general morphological models for complex fluids based on the completeness of these additive functionals).  Some drawbacks still  must be solved opening new research lines related to tractability of third and higher virial coefficients in a cluster expansion; the improvement of the accuracy of  threshold estimates based on Minkowski functionals; the computation of analytic expressions of Minkowski functionals beyond mean values of additive measures, and  the definition of these functionals on lattices  and their associated second--order moments
 (see  Sections 4 and 5 in \cite{Mecke}).

The additivity property of Minkowski functionals   becomes a limitation in systems exhibiting non-linear interactions,  non-linear phase behavior (e.g., phase separation in binary fluids or emulsions), fractal or self-similar structures (non-integer dimensions and scale-dependent structures), complex topologies (the presence of higher-order holes, handles, or cavities), overlapping particles  in granular and disordered systems, systems with non-additive energy contributions (e.g., curvature effects in membranes and vesicles involve non-linear energy terms), and complex percolation structures, where near the critical threshold,  clusters exhibit scale-free and non-linear growth.
 In these cases, additional tools, such as fractal analysis, non-additive geometric measures, or topological data analysis (TDA), may complement or replace the use of Minkowski functionals to better describe the system's behavior.

The present paper considers a more general family of  morphological measures allowing the dynamical analysis of spatial random structures  in physical systems, beyond the purely spatial analysis previously developed in the literature,  addressing the case where thresholds can change over time. The   asymptotic probabilistic  properties of this extended family of geometric functionals  allow, in particular, the characterization and prediction of future spatial configurations in random physical systems, providing basic quantities, like asymptotic  means and  variances, that play a crucial role  in  morphological analysis. In the context of percolation theory and Boolean model,  our approach introduces a more flexible stochastic modelling of  spatial patterns  incorporating  structural changes over time induced by time--varying thresholds.

It is worthwhile noticing  the  interest   in   Cosmology of the  extended asymptotic spatiotemporal analysis of Minkowski functionals addressed in this paper. Particularly, this analysis has special significance in Cosmic Microwave Background (CMB) Radiation Variation studies over time (see, e.g., \cite{Carones24}; \cite{Duque24}; \cite{Marinucci04}, among others). The CMB encodes information from the Early Universe in the intensity and polarisation of the light. The characterization of the statistical distribution of CMB  plays a crucial role in these studies.
 The analysis of morphological properties of  CMB spherical maps is based  on  completeness theorem,   the invariance assumption under translations and rotations, and on the additivity property of Minkowski functionals. Departure from  Gaussianity or deviations from isotropy assumption must be detected to validate statistical inference procedures, and to discriminate between competitive scenarios for Big Bang dynamics. The central role of Minkowski functionals in non--Gaussianity tests is well--known   (see, e.g., \cite{Marinucci04}).  In addition, Minkowski functionals  allow to address computational burden, the ease of masking or weighting data, and the analysis of deviations from different thresholds.  To go beyond the limited computation of  Minkowski functionals to CMB temperature and weak lensing, in \cite{Duque24}, Minkowski functionals are applied to  CMB polarisation data, introducing a new formalism that incorporates spin effects.
 CMB polarisation field is decomposed into two rotationally invariant fields, usually analyzed through their angular power spectra (see Section \ref{sec7} below). The  angular power spectrum  is  not sensitive to the possible deviations from Gaussianity or  departure from statistical isotropy. This fact motivates the extended application of Minkowski functionals to CMB polarisation data, exploiting their more efficient discrimination of Gaussian or isotropic deviations, with respect to the bispectrum, based on  three point correlation function, and trispectrum, based on four point correlation function  (see also  \cite{Carones24} for further insights in the analysis of CMB polarisation intensity maps from
Minkowski functionals). The asymptotic probability distribution results derived in the present paper, in particular, for spatiotemporal Minkowski functionals, provide new tools for the derivation of  a wider family of non--Gaussianity tests detecting  changes or deviations through time (see Section \ref{sec7} below).

From a statistical point of view,  the implementation of non--Gaussianity tests,  in the spatiotemporal context addressed in this paper, also requires the choice of a suitable parametric scenario for the family of spatiotemporal covariance functions, characterizing the two point correlation  structure of the underlying  spatiotemporal    isotropic random field. Note that variances of Minkowski  functionals also depend on the parameters of these covariance functions. In the context of
spatiotemporal Gaussian random fields (STGRFs),  the Gneiting class of  spatiotemporal covariance functions offers a flexible nonseparable   modelling framework
 (see  \cite{Gneiting02}).  An extended formulation  of Mat\'ern covariance function family to the spatiotemporal context can be achieved within this Gneiting class (see   \cite{Bevilacqua} for an unified framework in the purely spatial case).  See also Section \ref{examples} below, where additional examples of Gneiting class of nonseparable spatiotemporal  covariance functions are analyzed in the illustration of Theorem \ref{MFth1}, providing the asymptotic Gaussian distribution of spatiotemporal Minkowski functionals.

\subsection{State of the art}
Sojourn functionals were extensively analyzed since the nineties in
the context of weak--dependent random fields (see, e.g.,  \cite{Bulinski} and  \cite{Ivanov89}). Special
attention has been paid to the long--range dependent  random field case
(see  \cite{Ber};  \cite{Leonenko99};
\cite{LeonenkoOlenko14};  \cite{Makogin20};  \cite{MarinucciRV}, among others). Limit theorems for level
functionals of stationary Gaussian processes and fields constitute a major
topic in this literature (see, e.g.,  \cite{Auffinger13};  \cite{Wschebor09};
\cite{caba87};  \cite{Estrade16};
 \cite{Iribarren89};  \cite{KratzLeon}; \cite{KratzLeon2};  \cite{Marinucciv2};  \cite{Slud94}). These papers  contribute to the
characterization of topological and geometrical properties of random fields, from the analysis of morphological descriptors  like the Euler characteristic of an excursion set, the number of
up--crossings at level $u$ on a bounded closed cube of $\mathbb{R}^{d},$ the
probability of the maximum to be greater than a given threshold $u,$ or the
analysis of anisotropy based on the line integral with respect to the level
curve at any threshold  $u$ (see also     \cite{MullerD},  \cite{KratzV17} and  \cite{Pham13}).

 In the characterization of the asymptotic distribution of geometric functionals under LRD,   the limit results  derived are based on  reduction principle. This principle  was
first discovered by \cite{Taqqu} (see also \cite{Dobrushin};  \cite{Leonenkoetal17}; \cite{Leonenkoetal17b};  \cite{Taqqu2}).
Theorem 1 in \cite{LeonenkoRuiz-Medina2023}   provides a general reduction principle,  under increasing domain  asymptotics  in time, leading to the  limiting distributions of
properly normalised integral  functionals, when the underlying spatiotemporal Gaussian random field displays
 LRD in time.
   The restriction of this Gaussian field to the sphere  allows the application of these results in the context  of CMB analysis.  In this spherical context, interesting asymptotic results for sojourn functionals, under increasing domain asymptotics in time,   have been derived in \cite{MarinucciRV}, covering the cases where  the underlying spatiotemporal   random field displays Short--Range Dependence (SRD) and  LRD (see also
 \cite{Marinucciv2}). The special case of Hermite rank equal to two deserves special attention, since this case can be connected with Chi--squared statistics,   usually  arising  in non--Gaussianity tests in CMB analysis (see, e.g., \cite{Marinucci04}). See also
  \cite{Leonenkoetal17}.
 The results of this paper  extend the asymptotic analysis achieved in the above cited references  to the case of  increasing domain  asymptotics  in time and space, incorporating moving thresholds   (see
Section \ref{contribution} below for more technical details).

\subsection{Our contribution}
\label{contribution}

Up to our knowledge,  the problem of determining, via reduction theorems, the asymptotic probability distribution  of integral functionals of nonlinear transformations of LRD  STGRFs, under   increasing domain asymptotics in space and time, has not been addressed in   the current literature. The extended analysis   to the case of time--varying
subordination, via a nonlinear transformation depending on the size of the temporal
domain  constitutes the main contribution of the present paper. First Minkowski functionals involving moving
levels subordinated to  STGRFs  arise as interesting particular case of this extended analysis . Namely, we restrict our attention to the case where time domain is a temporal continuous  interval $[0,T],$ and $\Lambda
(T)=T^{\gamma },$ $\gamma \geq 0,$ defines the scale factor
of the homothetic transformation  of a convex compact set $K\subset \mathbb{R%
}^{d},$ with center at point $\mathbf{0}\in K.$   The case $\gamma = 0$  was addressed  in  \cite{LeonenkoRuiz-Medina2023} (see also  \cite{Marinucciv2}, where alternative conditions were formulated in  a Hilbert space framework, under a parametric modelling of spherical scale dependent LRD). When $\gamma =1,$ our Theorem \ref{th3} extends Theorem 1 in
  \cite{Leonenkoetal17} to the spatiotemporal context (see also    \cite{Leonenkoetal17b}).

In the second part of the paper,  Theorem \ref{th2} provides  a general reduction
principle under $T$-varying nonlinear transformations of Gaussian random
fields, with, as before, $T$ denoting the size of the time interval. Additionally to the assumed conditions on LRD in the first part of the paper, in the derivation of
Theorem \ref{th3b}, the divergence rate of the moving threshold parameter must be controlled at
the logarithm scale by the increasing of the size $T$ of the temporal
interval (see \textbf{Condition 6} below).  In that sense, the   methodological approach adopted in the  proof  of this result  is not standard.  This reduction principle  characterizes the limit Gaussian distribution of Minkowski functionals
involving moving levels.
The application of the above results to the context of sojourn measures of  spatiotemporal Gaussian random fields restricted to the  sphere is then contemplated,
motivated  by the current literature on CMB radiation variation analysis.

The outline of the paper is now introduced.  Some preliminary results on geometrical
probabilities are first reviewed in Section \ref{s2}. Section \ref{sec3}
derives the conditions for a general reduction principle under increasing
domain asymptotics in time and in space (see Theorem \ref{th3}). In Section
\ref{exsec}, Theorem \ref{MFth1} derives the asymptotic Gaussian distribution of spatiotemporal
Minkowski functionals by applying reduction principle provided in  Theorem \ref{th3}.
Some examples are also analyzed, where subordination to
 STGRFs with nonseparable covariance function is considered.
The case where  geometric  integral  functionals  are computed from a $T$--dependent nonlinear
transformations of a Gaussian random field is studied in Section \ref{ef}
(see Theorem \ref{th2}). In Section \ref{s5}, the asymptotic distribution of  suitable
normalized first Minkowski functionals involving moving levels, under
increasing domain asymptotics in space--time, is derived in Theorem \ref{th3b}.
 The obtained  results  are applied in  Section \ref{sec7} to the case  of $\gamma =0,$
 when the underlying Gaussian spatiotemporal random field  is  the restriction to the sphere of a  stationary in time, and homogeneous and isotropic in space, STGRF.

\section{Preliminaries}
\label{s2} Under  spatial   isotropy,
variance components in the chaotic expansion of spatiotemporal  geometric  integral  functionals can be  computed from geometrical probabilities
 (see, e.g.,  \cite{Lord54};  \cite{Ivanov89}). Some extended results  can also  be found  in
\cite{Aharonyan20},  and the references
therein.  We now summarize them and introduce the corresponding notation.

Let $\nu_{d}(\cdot)$ be the Lebesgue measure in $\mathbb{R}^{d},$ $%
d\geq 2,$ and $K$ be a convex body in $\mathbb{R}^{d},$ i.e., a compact
convex set with non--empty interior with center at the point $0\in K.$ Let $%
\mathcal{D}(K)=\max\left\{ \|x-y\|,\ x, y\in K\right\}$ be the diameter of $%
K.$ Let also $\nu_{d}(K)=|K|$ be the volume of $K,$ and $\nu_{d-1}(\delta
K)=U_{d-1}(K)$ be the surface area of $K,$ where $\delta K$ is the boundary
of $K.$  \ Note that for  $K=\mathcal{B}(1)=\{x\in \mathbb{R}^{d};\
\|x\|\leq 1\}$ is the unit ball, and $\delta K=S_{d-1}(1)=\{x\in \mathbb{R}%
^{d};\ \|x\|=1\}$ is the unit sphere, then $\mathcal{D}(\mathcal{B}(1))=2,$
and  $|\mathcal{B}(1)|=\frac{\pi^{d/2}}{\Gamma \left(\frac{d}{2}+1\right)},\
U_{d-1}(\mathcal{B}(1))=|S_{d-1}|=\frac{2\pi^{\frac{d}{2}}}{\Gamma \left(%
\frac{d}{2}\right)}.$

Let $\Lambda (T)K$ be a homothetic transformation of body $K$ with
center $0\in K,$ and coefficient or scaling factor $\Lambda (T)>0.$
We assume that  $\Lambda (T)=T^{\gamma },$ $\gamma \geq
0 $. Then, $|\Lambda (T)K|=|K|T^{\gamma d}.$ Following the approach
presented in  \cite{Lord54} (see also  \cite{Ivanov89}%
), the probability density $\psi _{\Lambda (T),\mathcal{B}(1)}(z)$ of the
random variable $Z=\rho_{\mathcal{B}(1)}=\|P_{1}-P_{2}\|,$ with $P_{1}$ and $%
P_{2}$ being two independent random points with uniform distribution in $%
\mathcal{B}(1),$ is given by, for $0\leq z\leq 2T^{\gamma},$

\begin{eqnarray}
\psi_{\Lambda (T),\mathcal{B}(1)}(z)&=&\frac{d}{[\Lambda (T)]^{d}}%
z^{d-1}I_{1-\left(\frac{z}{2\Lambda (T)}\right)^{2}} \left(\frac{d+1}{2},%
\frac{1}{2}\right),  \label{b1bb}
\end{eqnarray}
\noindent in terms of the incomplete beta--function

\begin{equation}
I_{\mu}(p,q)=\frac{\Gamma (p+q)}{\Gamma (p)\Gamma (q)}\int_{0}^{%
\mu}t^{p-1}(1-t)^{q-1}dt,\quad \mu\in (0,1].  \label{ibf}
\end{equation}

In formula (2.6) in  \cite{Aharonyan20}, an extended
version of the probability density of the distance between two independent
uniformly distributed points in  a convex body $K$ in $\mathbb{R}^{d}$ is derived,  applying an alternative methodology to  \cite{Lord54} for the
case of hyperspheres. Specifically, let $\mathcal{J}$ be the space of all straight lines in $\mathbb{R}^{d},$  and $d\gamma $ is an element of a locally finite measure in the space
$\mathcal{J},$ which is invariant, with respect to the group $M$ of all Euclidean motions in the space $\mathbb{R}^{d}$  (the uniform measure on $\mathcal{J}$). Let  also $F_{K}(v)$ be the chord length distribution function of body $K,$  defined as
$$F_{K}(v)=\frac{2(d-1)}{|S_{d-2}|}\int_{|\chi(\gamma )|\leq v }d\gamma,$$
\noindent where $\chi(\gamma )=\gamma \cap K$ is a chord in $K.$ Then, for $0\leq z\leq \mathcal{D}(K),$
\begin{eqnarray}
\psi_{\rho_{K}}(z)=\frac{1}{|K|^{2}}\left[%
z^{d-1}|S_{d-1}||K|-z^{d-1}|S_{d-2}|\frac{U_{d-1}(K)}{d-1}%
\int_{0}^{z}\left(1-F_{K}(v)\right)dv\right].  \label{eqfdfd}
\end{eqnarray}

 From (\ref{eqfdfd}), we also have, for $0\leq z\leq
\mathcal{D}(\Lambda (T)K),$

\begin{eqnarray}
&&\psi_{\Lambda (T),K}(z)=\psi_{\rho _{\Lambda (T)K}}(z)=\frac{1}{%
[\left(|K|[\Lambda (T)]^{d}\right)]^{2}}\left[\left(|K|[\Lambda
(T)]^{d}\right)z^{d-1}|S_{d-1}|\right.  \notag \\
&& \hspace*{3.5cm}\left.-|S_{d-2}|z^{d-1} \frac{U_{d-1}(\Lambda (T)K)}{d-1}%
\int_{0}^{z}\left(1-F_{\Lambda (T)K}(v)\right)dv\right].  \label{pdd}
\end{eqnarray}
\section{Reduction theorems for spatiotemporal random fields with LRD}
\label{sec3}
We consider the spatiotemporal random field
 $Z:\left( \Omega \times \mathbb{R}^{d}\times \mathbb{R}\right) \longrightarrow \mathbb{R},$
with  $(\Omega ,\mathcal{A},P)$  denoting  the basic probability space.

 \medskip

\noindent \textbf{Condition 1}. Assume that  $Z$ is  a measurable mean--square
continuous homogeneous and isotropic in space and stationary in time
Gaussian random field with
\begin{eqnarray}
\mathbb{E}[Z(x,t)]&=& 0, \quad \mathbb{E}[Z^{2}(x,t)]=1 \nonumber\\
\widetilde{C}(\|x-y\|,|t-s|)&=&\mathbb{E}\left[Z(x,t)Z(y,s)\right]\geq 0.\nonumber
\end{eqnarray}

\medskip

Under \textbf{Condition 1}, one can write
\begin{eqnarray}
C(z,\tau) &=&\widetilde{C}(\|x-y\|,|t-s|),\quad z=\|x-y\|\geq 0, \quad
\tau=|t-s|\geq 0,\nonumber\end{eqnarray}

\noindent where $C$ denotes the covariance function as a function of the arguments  $z$ and $\tau,$ respectively representing the norm of the spatial argument, and the absolute value of the  time argument.  While  $\widetilde{C}$ means that we are considering the values of the covariance function depending on the input arguments $x-y\in \mathbb{R}^{d}$ and $t-s\in \mathbb{R}.$

\begin{remark}
The non-negativeness condition of the covariance function has been usually assumed in the literature of LRD stationary  Gaussian processes,  and isotropic  random fields during nineties, as standard assumption in  reduction theorems, leading to central and non--central limit results for integral functionals  of nonlinear transformations of Gaussian random fields.  For an alternative approach in the derivation of limit results for weighted nonlinear transformations of  LRD Gaussian stationary processes  we refer the reader to \cite{Savich13}.  See also
\cite{Mainia24}, where  the application of the Malliavin–Stein method and Fourier analysis techniques in the derivation of spectral limit theorems is considered. The general setting of product of spatial domains of different dimensions  requiring  different scaling factors   is also analyzed in \cite{LeonenkoMP2024}.
\end{remark}

Let $Z\sim \mathcal{N}(0,1)$ be a standard Gaussian random variable
with density
\begin{equation*}
\phi (w)=\frac{1}{\sqrt{2\pi }}\exp\left(-\frac{w^{2}}{2}\right),\ w\in%
\mathbb{R},\ \Phi (u)=\int_{-\infty}^{u}\phi (w)dw.
\end{equation*}

Let $G(z)$ be a real--valued Borel function satisfying

\medskip

\noindent \textbf{Condition 2}. $\mathbb{E}G^{2}(Z(x,t))<\infty .$

\medskip

Under \textbf{Condition 2},  $G$  admits a chaotic  expansion  in terms of the normalized
Hermite polynomials in the Hilbert space $L_{2}(\mathbb{R},\phi (u)du)$ of square integrable functions with respect to the standard Gaussian measure. This expansion is given by
\begin{equation*}
G(z)=\sum_{n=0}^{\infty }\frac{\mathcal{J}_{n}}{n!}H_{n}(z),\ \mathcal{J}%
_{n}=\int_{\mathbb{R}}G(z)H_{n}(z)\phi (z)dz,
\end{equation*}%
\noindent where the Hermite polynomial of order $n,$ denoted by $H_{n}(z),$
is defined by the equation \linebreak $\frac{d^{n}}{dz^{n}}\phi (z)=(-1)^{n}H_{n}(z)\phi (z).$
Note that $H_{0}=1,$ $H_{1}(z)=z,$ and $H_{2}(z)=z^{2}-1,\dots $

\medskip

Following \cite{Taqqu}, we will introduce the following condition:

\medskip

\noindent \textbf{Condition 3}. The Hermite rank of the function $G$
is $m\geq 1.$ Hence, for $m=1,$ $\mathcal{J}_{1}\neq 0,$ or for $m\geq 2,$ $%
\mathcal{J}_{1}=\cdots=\mathcal{J}_{m-1}=0,$ $\mathcal{J}_{m}\neq 0$ (see also  \cite{Taqqu2}).

Under \textbf{Condition 3},  the following spatiotemporal integral functional of
$Z(x,t)$ is defined:

\begin{eqnarray}
A_{T}&=&\int_{0}^{T}\int_{\Lambda (T)K}G(Z(x,t))dxdt=\mathcal{J}_{0}+\sum_{n\geq m}%
\frac{\mathcal{J}_{n}}{n!}\xi_{n,T}  \nonumber \\
E[A_{T}]&=&\int_{0}^{T}\int_{\Lambda (T)K}\mathcal{J}_{0}dxdt=\mathcal{J}%
_{0}|T||K|[\Lambda (T)]^{d}  \nonumber \\
\xi_{n,T}&=& \int_{0}^{T}\int_{\Lambda (T)K}H_{n}(Z(x,t))dxdt  \notag \\
\mathbb{E}[\xi_{n,T}] &=& 0,\ \mathbb{E}[\xi_{n,T}\xi_{l,T}]=0,\ n\neq l,\
n,l\geq m,  \label{cond}
\end{eqnarray}
\noindent where, as given in \textbf{Condition 3},  $m$ denotes the Hermite rank of function $G,$ and the integrals are interpreted in the mean square sense.

Thus, for $n\geq m,$
\begin{eqnarray}  \label{nevar2bb}
&&\sigma _{n,K}^{2}(T)=\mbox{Var}(\xi _{n,T})  \notag \\
&=& 2n!T\int_{0}^{T}\left(1-\frac{\tau}{T}\right)\int_{K\Lambda (T)\times
K\Lambda (T)} \widetilde{C}^{n}(\|x-y\|,\tau)d\tau d%
xdy  \notag \\
&=&2n!T|K\Lambda (T)|^{2}\int_{0}^{T}\left(1-\frac{\tau}{T}\right)\mathbb{E}%
\left[\widetilde{C}^{n}\left(\|P_{1}-P_{2}\|,\tau\right) \right]d\tau  \notag
\\
&&= 2n!T|K|^{2}T^{2\gamma d}\int_{0}^{T}\left( 1-\frac{\tau }{T}\right)
\int_{0}^{\mathcal{D}(K\Lambda (T))}\psi _{\Lambda (T),K}(z)C^{n}(z,\tau
)dzd\tau .  \notag \\
\end{eqnarray}
\noindent Here, $\psi _{\Lambda (T),K}(z),$ denotes, as before,  the probability
density of the random variable \linebreak $\rho_{\Lambda (T)K}=\|P_{1}-P_{2}\|,$ where $P_{1}$
and $P_{2}$ are two independent random points with uniform distribution in $%
\Lambda (T)K.$ In particular, from equation (\ref{b1bb}), for $K=\mathcal{B}(1),$ for $0\leq z\leq 2T^{\gamma
},$

\begin{eqnarray}
\psi_{T^{\gamma
},\mathcal{B}(1)}(z)&=&\frac{d}{T^{\gamma d}}%
z^{d-1}I_{1-\left(\frac{z}{2T^{\gamma }}\right)^{2}} \left(\frac{d+1}{2},%
\frac{1}{2}\right),  \label{b1}
\end{eqnarray}
\noindent in terms of the incomplete beta function (\ref{ibf}).
Thus,

\begin{eqnarray}
\sigma _{m,\mathcal{B}(1)}^{2}(T) &=& \mbox{Var}\left(
\int_{0}^{T}\int_{\Lambda (T)\mathcal{B}(1)} H_{m}\left(
Z(x,t)\right)\right)dxdt \nonumber\\ &=& 2m! T|\mathcal{B}(1)|^{2}[\Lambda (T)]^{2d}\int_{0}^{T}\left( 1-\frac{%
\tau }{T}\right) \int_{0}^{\mathcal{D}(K\Lambda (T))}C^{m}(z,\tau )  \nonumber \\
& &\hspace*{1cm}\times \left[\frac{d}{[\Lambda (T)]^{d}}z^{d-1}I_{1-\left(%
\frac{z}{2\Lambda (T)}\right)^{2}} \left(\frac{d+1}{2},\frac{1}{2}\right)%
\right]dzd\tau  \nonumber \end{eqnarray}

\begin{eqnarray}
 &=& 2m!|\mathcal{B}(1)|^{2}dT^{\gamma d+1}\int_{0}^{T}\left( 1-\frac{\tau }{T%
}\right)  \nonumber\\
&&\hspace*{2cm}\times \int_{0}^{2T^{\gamma }}z^{d-1}C^{m}(z,\tau
)I_{1-\left( \frac{z}{2T^{\gamma }}\right) ^{2}}\left( \frac{d+1}{2},\frac{1%
}{2}\right) dzd\tau  \nonumber \\
&=& 8m!\frac{\pi^{d}}{d\Gamma^{2}\left( \frac{d}{2}\right)}T^{\gamma
d+1}\int_{0}^{T}\left( 1-\frac{\tau }{T}\right) \int_{0}^{2T^{\gamma
}}z^{d-1}C^{m}(z,\tau )  \nonumber \\
&& \hspace*{3.5cm}\times I_{1-\left( \frac{z}{2T^{\gamma }}\right)
^{2}}\left( \frac{d+1}{2},\frac{1}{2}\right) dzd\tau.  \label{eqbvc4}
\end{eqnarray}

\begin{remark}
\label{remlord} In the following, we will apply Lord (1954) results on the
derivation of the probability density of the distance between two
independent uniform points in hypersheres. Consider the positive constants $%
F_{1}$ and $F_{2}$ respectively defining the supremum and infimum of the
radius of the balls such that the following inclusions hold
\begin{equation}
F_{1}\mathcal{B}\left(1\right)\subseteq K\subseteq F_{2}\mathcal{B}\left(
1\right).\label{cb}
\end{equation}

From equations (4)--(6) in Lord (1954),
\begin{equation}
C_{1}\psi_{T^{\gamma},\mathcal{B}\left(S_{1}\right)}\leq
\psi_{T^{\gamma},K}\leq C_{2}\psi_{T^{\gamma},\mathcal{B}\left(
S_{2}\right)},\ 0\leq C_{1}\leq C_{2}.  \label{ineqlord}
\end{equation}
\noindent Here, $\psi _{T^{\gamma},K}(z)$ denotes the probability density of
the random variable $\rho_{KT^{\gamma }}=\|P_{1}-P_{2}\|,$ with $P_{1}$
and $P_{2}$ being two independent random points with uniform distribution in
$KT^{\gamma }.$ As before, $K$ denotes a compact convex set with non--empty
interior, and with center at the point $0\in K.$ In particular, $\mathcal{B}%
\left(S_{i}\right)$ denotes the ball with center $0\in \mathcal{B}%
\left(S_{i}\right),$ and radius $S_{i},$ $i=1,2,$ according to constants $F_{i},$ $i=1,2,$ in equation
(\ref{cb}).
\end{remark}

The following additional condition will be assumed.

\noindent \textbf{Condition 4}.

\begin{itemize}
\item[(i)] $C(z,\tau )\rightarrow 0,$ if $\max \{z,\tau
\}\rightarrow \infty .$

\item[(ii)] For some fixed $m\in \{1,2,3\dots \},$ there exist $\delta _{1}\in (0,1)$ and $\delta _{2}\in (0,1)$ such that for
$\Lambda (T)=T^{\gamma },$ for certain $\gamma
\geq 0,$
\begin{eqnarray}
&& \lim_{T\rightarrow \infty }\frac{\sigma _{m,K}^{2}(T)}{T^{1+\delta _{1}}T^{\gamma d(1+\delta _{2})}    }=\infty.  \nonumber
\end{eqnarray}

\end{itemize}

\medskip

From Remark \ref{remlord} (see also equations (\ref{nevar2bb}), (\ref{eqbvc4}) and (\ref{ineqlord})), it is straightforward that \textbf{Condition 4(ii)} holds if there exist $\delta _{1}\in
(0,1)$ and $\delta _{2}\in (0,1)$ such that as $T\to \infty,$

\begin{eqnarray}
&&\frac{1}{T^{\delta _{1}+\gamma d\delta _{2}}}\int_{0}^{T} \left( 1-\frac{\tau }{T}\right) \int_{0}^{2T^{\gamma }}%
z^{d-1}C^{m}(z,\tau )
I_{1-\left( \frac{z}{2T^{\gamma }}\right) ^{2}}\left( \frac{d+1}{2},\frac{1}{%
2}\right) dzd\tau\to \infty .  \nonumber\\
\label{ineqversus}
\end{eqnarray}

In the subsequent development, \textbf{Condition 4(ii)} will be  verified in terms of
equation (\ref{ineqversus}).

\begin{remark}
\label{inter}
\textbf{Condition 4}  and equation (\ref{ineqversus}) mean that the spatiotemporal Gaussian random field $Z$ displays LRD in space and time. The general reduction principle provided in Theorem \ref{th3} below essentially follows from this condition, for any function $G\in L_{2}(\mathbb{R},\phi(u)du),$ under the Gaussian scenario introduced in \textbf{Condition 1} (see also equation (\ref{eqc4ii}) in the proof of this result). This fact is illustrated in Sections \ref{separablecov} and \ref{examples}. Specifically, the case of separable covariance functions in space and time is considered in  Section \ref{separablecov}.
In this  case,  \textbf{Condition 4(ii)} holds for $0<A<1/m,$ and $0<\widetilde{\alpha }<d/m$    (see equations  (\ref{lrdtime}) and  (\ref{LRDspace})), which corresponds to the case of  LRD  in time and space.  The nonseparable covariance function case is illustrated in Section  \ref{examples} within the
Gneiting covariance function class (see equations (\ref{gneit})--(\ref{fcf})). The sufficient conditions considered,
 $0<\delta_{2}<1-\frac{\alpha \beta}{\gamma}<1,$ $
0<\delta_{1}<1-2\widetilde{\gamma } (\gamma -2\alpha \beta )<1,$ in the first example,  and
$0<\delta_{2}<1-\frac{\alpha \beta}{\gamma}<1,$ $
0<\delta_{1}<1-2\nu\widetilde{\gamma } (\gamma -2\alpha \beta )<1,$  in the second  example, reflect LRD in time and space,  involving  space--time interaction  in the restrictions on the LRD parameters $\alpha , \beta $ (time), $\widetilde{\gamma },\nu$ (space), as well as on the shape parameter $\gamma $ characterizing the scaling factor of the homothetic transformation  of   $K\subset \mathbb{R%
}^{d}.$ Note that, in both cases separable and nonseparable covariance function cases, the  LRD parameter are  involved in  the variance of
Minkowski functionals.

\end{remark}

%%%%%%%
\subsection{Separable covariance functions}
\label{separablecov}
\label{scf}
For the case of separable spatiotemporal covariance functions in the unit ball,
\begin{equation*}
C(z,\tau )=C_{Space}(z)C_{Time}(\tau ),
\end{equation*}%
 we have
\begin{equation}
\sigma _{m,\mathcal{B}(1)}^{2}(T)=m!b_{1m}(T)b_{2m}(T),
\label{sigmam}
\end{equation}%
\noindent where
\begin{eqnarray}
&& b_{1m}(T)=
2T\int_{0}^{T}\left( 1-\frac{\tau }{T}\right) C_{Time}^{m}(\tau
)d\tau  \nonumber\\
&& b_{2m}(T)=|\mathcal{B}(1)|^{2}dT^{\gamma d}\int_{0}^{2T^{\gamma
}}z^{d-1}C_{Space}^{m}(z)I_{1-\left( \frac{z}{2T^{\gamma }}\right)
^{2}}\left( \frac{d+1}{2},\frac{1}{2}\right) dz.\nonumber
\end{eqnarray}

Then, in the weak dependent case in time, i.e., when the temporal covariance function is absolutely integrable,  as T  $\rightarrow \infty ,$ we obtain
\begin{eqnarray}
&& b_{1m}(T)=2L_{1}T(1+o(1))\nonumber\\
&&
L_{1}=\int_{0}^{\infty }C_{Time}^{m}(\tau )d\tau <\infty ,\quad
\int_{0}^{\infty }C_{Time}^{m}(\tau )d\tau \neq 0.
\label{wdc}
\end{eqnarray}%

Suppose that for some bounded slowly varying functions at
infinity $\mathcal{L}_{1}(\tau ):$
\begin{equation}
C_{Time}(\tau )=\frac{\mathcal{L}_{1}(\tau )}{\tau ^{A}},\quad A>0.
\label{lrdtime}
\end{equation}%
\noindent
Then,
\begin{eqnarray}
&&b_{1m}(T)
=2T\int_{0}^{T}\left(1-\frac{\tau}{T}\right)\left[\frac{\mathcal{L}_{1}(\tau )}{\tau^{A}}\right]^{m}d\tau.\nonumber
\end{eqnarray}

Under covariance model (\ref{lrdtime}),  equation (\ref{wdc}) holds for $A>1/m,$ corresponding to the weak--dependent case in time.
For $A=\frac{1}{m},$
$b_{1m}(T)=
2T\log (T)\mathcal{L}_{1}^{m}(T)(1+o(1)).$  While for  $0<A<\frac{1}{m},$ applying  the change of variable $x=\frac{\tau}{T},$
\begin{eqnarray}
b_{1m}(T)&=& 2T^{2-mA}\mathcal{L}_{1}^{m}(T)
\int_{0}^{1}\frac{\mathcal{L}_{1}(xT )}{\mathcal{L}_{1}^{m}(T)}\frac{1}{x^{Am}}(1-x)dx\nonumber\\
&=&
2L_{2}T^{2-mA}\mathcal{L}_{1}^{m}(T)(1+o(1))  \nonumber\\
L_{2}&=&\left[ \int_{0}^{1}(1-\tau )\tau ^{-Am}d\tau \right] =\left[(1-mA)(2-Am)\right]^{-1}.
\label{eqrev1}
\end{eqnarray}

Similarly, as  $T\rightarrow \infty ,$ under weak--dependence in space, i.e., when the spatial covariance function is absolutely integrable over the spatial domain,  $$b_{2m}(T)=L_{3}T^{d\gamma }(1+o(1)),$$
\noindent where \begin{eqnarray}L_{3}&=&|\mathcal{B}(1)|^{2}d\int_{0}^{\infty
}z^{d-1}C_{Space}^{m}(z)dz<\infty \nonumber\\ &&\int_{0}^{\infty
}z^{d-1}C_{Space}^{m}(z)dz\neq 0.\nonumber\end{eqnarray}

Furthermore,   if for some bounded slowly varying functions at
infinity   $\mathcal{L}(z):$
\begin{equation}
C_{Space}(z)=\frac{\mathcal{L}(z)}{z ^{\widetilde{\alpha } }},\quad
\widetilde{\alpha }>0,
\label{LRDspace}
\end{equation}
\noindent then, for $\widetilde{\alpha }=\frac{d}{m},$
\begin{eqnarray}
&&b_{2m}(T)=|\mathcal{B}(1)|^{2}dT^{\gamma d}\int_{0}^{2T^{\gamma
}}z^{d-1}\left[\frac{\mathcal{L}(z)}{z ^{\widetilde{\alpha } }}\right]^{m}I_{1-\left( \frac{z}{2T^{\gamma }}\right)
^{2}}\left( \frac{d+1}{2},\frac{1}{2}\right) dz\nonumber\\
&&=
L_{4}T^{\gamma d}\log (T^{\gamma })\mathcal{L}^{m}(T^{\gamma
})(1+o(1)), \nonumber\end{eqnarray}
\noindent
with $L_{4}=\frac{4\pi^{d}}{d\Gamma^{2}(d/2)}.$ For $0<\widetilde{\alpha } <\frac{d}{m},$ considering
the change of variable $u=\frac{z}{T^{\gamma }}$
\begin{eqnarray}
&&b_{2m}(T)=|\mathcal{B}(1)|^{2}dT^{\gamma d}\int_{0}^{2T^{\gamma
}}z^{d-1}\left[\frac{\mathcal{L}(z)}{z ^{\widetilde{\alpha } }}\right]^{m}I_{1-\left( \frac{z}{2T^{\gamma }}\right)
^{2}}\left( \frac{d+1}{2},\frac{1}{2}\right) dz\nonumber\\
&&=|\mathcal{B}(1)|^{2}dT^{\gamma (2d-m\widetilde{\alpha })}\mathcal{L}^{m}(T^{\gamma
})\int_{0}^{2}u^{d-1}\left(\frac{\mathcal{L}^{m}(T^{\gamma}u)}{\mathcal{L}^{m}(T^{\gamma
})}-1+1\right)I_{1-\left( \frac{u}{2}\right)
^{2}}\left( \frac{d+1}{2},\frac{1}{2}\right) du\nonumber\\
&&=|\mathcal{B}(1)|^{2}dT^{\gamma (2d-m\widetilde{\alpha })}\mathcal{L}^{m}(T^{\gamma
})\left[ \int_{0}^{2}u^{d-1}I_{1-\left( \frac{u}{2}\right)
^{2}}\left( \frac{d+1}{2},\frac{1}{2}\right) du\right.\nonumber\\
&&\left.\hspace*{2cm}+\int_{0}^{2}u^{d-1}\left(\frac{\mathcal{L}^{m}(T^{\gamma}u)}{\mathcal{L}^{m}(T^{\gamma
})}-1\right)I_{1-\left( \frac{u}{2}\right)
^{2}}\left( \frac{d+1}{2},\frac{1}{2}\right) du\right]\nonumber\\
&&=L_{5}T^{\gamma (2d-m\widetilde{\alpha } )}\mathcal{L}%
^{m}(T^{\gamma })(1+o(1)), \label{eqrev2}
\end{eqnarray}
\noindent
 \noindent where
$$L_{5}=\frac{2^{d-m\widetilde{\alpha } +1}\pi ^{d-\frac{1}{2}}\Gamma \left(
\frac{d-m\widetilde{\alpha } +1}{2}\right) }{(d-m\widetilde{\alpha } )\Gamma
\left( \frac{d}{2}\right) \Gamma \left( \frac{2d-m\widetilde{\alpha }+2}{2}%
\right) },$$ \noindent  since
\begin{eqnarray}
&&\int_{0}^{2}u^{d-1}I_{1-\left( \frac{u}{2}\right)
^{2}}\left( \frac{d+1}{2},\frac{1}{2}\right) du = 2^{d}\int_{0}^{1}w^{d-1}I_{1-w^{2}}\left( \frac{d+1}{2},\frac{1}{2}\right)dw\nonumber\\
&&=\frac{2^{d}}{B\left(\frac{d+1}{2},\frac{1}{2}\right)}\int_{0}^{1}\int_{0}^{1-w^{2}}w^{d-1}t^{\frac{d-1}{2}}(1-t)^{-1/2}dtdw\nonumber\\
&&=\frac{2^{d}}{B\left(\frac{d+1}{2},\frac{1}{2}\right)}\int_{0}^{1}t^{\frac{d-1}{2}}(1-t)^{-1/2}\left[\int_{0}^{\sqrt{1-t}}w^{d-1}dw\right]dt
\nonumber\\
&&=\frac{2^{d}}{dB\left(\frac{d+1}{2},\frac{1}{2}\right)}\int_{0}^{1}(1-t)^{(d-1)/2}t^{\frac{d-1}{2}}dt\nonumber\\
&&=\frac{2^{d}B\left(\frac{d+1}{2},\frac{d+1}{2}\right)}{dB\left(\frac{d+1}{2},\frac{1}{2}\right)},\nonumber
\end{eqnarray}
\noindent and \begin{equation}\lim_{T\to \infty}\int_{0}^{2}u^{d-1}\left(\frac{\mathcal{L}^{m}(T^{\gamma}u)}{\mathcal{L}^{m}(T^{\gamma
})}-1\right)I_{1-\left( \frac{u}{2}\right)
^{2}}\left( \frac{d+1}{2},\frac{1}{2}\right) du=0\label{eqrev2b}\end{equation}
 \noindent (see also Lemma 2.1.3 of  \cite{Ivanov89}).

Summarizing,  for the introduced separable covariance function family,  when $K=
\mathcal{B}(1),$ \textbf{Condition 4(ii)} holds in the following
cases:

From equation (\ref{sigmam}), when  $\gamma =0,$ and  $0<A<\frac{1}{m},$
\begin{equation}
\sigma _{m,\mathcal{B}(1)}^{2}(T) =m!L_{2}T^{2-mA}\mathcal{L}
_{1}^{m}(T)(1+o(1)),\label{eqrev3}\end{equation}
\noindent hence, \begin{eqnarray}
&& \lim_{T\rightarrow \infty }\frac{\sigma _{m,\mathcal{B}(1)}^{2}(T)}{T^{1+\delta _{1}}}   =\infty,  \label{eqrev3b}
\end{eqnarray}
\noindent  for  $0<\delta _{1}<1-mA<1,$   in
the case of LRD in time. This case has been considered by
\cite{LeonenkoRuiz-Medina2023}. While for $\gamma >0,$ $0<A<\frac{1}{m},$  and $0<\widetilde{\alpha } <\frac{d}{m},$
\begin{eqnarray} &&\sigma _{m,\mathcal{B}(1)}^{2}(T) =m!L_{2}L_{5}T^{2-mA}\mathcal{L}
_{1}^{m}(T)T^{\gamma (2d-m\widetilde{\alpha } )}\mathcal{L}^{m}(T^{\gamma
})(1+o(1)),\label{eqrev4}\end{eqnarray}
\noindent and \begin{eqnarray}
&& \lim_{T\rightarrow \infty }\frac{\sigma _{m,\mathcal{B}(1)}^{2}(T)}{T^{1+\delta _{1}}T^{\gamma d(1+\delta _{2})}    }=\infty,  \nonumber
\end{eqnarray}
\noindent  for $0<\delta _{1}<1-mA<1,$
and $0<\delta _{2}<1-(m\widetilde{\alpha } )/d<1,$
in the case of  LRD in time and space.  For the remaining cases, \begin{equation}\lim_{T\rightarrow \infty }\frac{\sigma _{m,\mathcal{B}(1)}^{2}(T)}{T^{1+\delta _{1}}T^{\gamma d(1+\delta _{2})}    }=0,\label{eqrev5}\end{equation}
\noindent for any $\delta_{i}>0,$ $i=1,2,$ and  \textbf{Condition 4(ii)}  does not hold.

\begin{theorem}
\label{th3} Under  \textbf{Conditions 1,2,3 and 4}, the random
variables
\begin{eqnarray}
Y_{T} &=&\frac{\left[ A_{T}-\mathbb{E}[A_{T}]\right] }{|\mathcal{J}%
_{m}|\sigma _{m,K}(T)/m!}  \label{funct1} \\
&&  \notag \\
&&\hspace*{-3.75cm}\mbox{and}  \notag \\
&&  \notag \\
&&Y_{m,T}=\frac{\mbox{sgn}(\mathcal{J}_{m})\int_{0}^{T}\int_{T^{\gamma
}K}H_{m}(Z(x,t))dxdt}{\sigma _{m,K}(T)}  \label{funct2}
\end{eqnarray}%
\noindent have the same limiting distributions as $T\rightarrow \infty $ (if
one of them  exists).
\end{theorem}

\noindent \textbf{Proof.} We consider the decomposition
$A_{T}-\mathbb{E}[A_{T}]=S_{1,T}+S_{2,T},$  where
\begin{eqnarray}
S_{1,T}&=&\frac{\mathcal{J}_{m}}{m!}\xi_{m,T} ,\quad
S_{2,T}= \sum_{n=m+1}^{\infty}\frac{\mathcal{J}_{n}}{n!}\xi_{n,T},\
\sum_{n=m}^{\infty}\frac{\mathcal{J}^{2}_{n}}{n!}<\infty.  \label{eqs1s2}
\end{eqnarray}

From (\ref{cond}),
\begin{equation}
\mbox{Var}(A_{T})=\mbox{Var}[S_{1,T}]+\mbox{Var}[S_{2,T}],
\label{eqvarfunct}
\end{equation}%
\noindent and we will prove that $\mbox{Var}[S_{2,T}]/\sigma _{m,K\Lambda
(T)}^{2}(T)\rightarrow 0,$ $T\rightarrow \infty .$

From \textbf{Condition 4(i)},
\begin{eqnarray}
&&\sup_{(z,\tau )\in B_{T}^{\delta _{1},\delta _{2}}}C(z,\tau )\rightarrow
0,\quad T\rightarrow \infty  \nonumber\\
&&B_{T}^{\delta _{1},\delta _{2}}=\{(z,\tau );\ \tau \geq T^{\delta _{1}}\ %
\mbox{or}\ z\geq T^{\gamma \delta _{2}}\}.  \label{cz}
\end{eqnarray}

Furthermore, for the set
\begin{equation*}
\overline{B}_{T}^{\delta _{1},\delta _{2}}=\left\{ (z,\tau);\ 0\leq \tau
\leq T^{\delta _{1}},\ 0\leq z\leq T^{\gamma \delta _{2}}\right\} ,
\end{equation*}%
\noindent one can use the estimate $C^{m+1}(z,\tau )\leq 1.$ Let us
consider $\delta _{1}\in (0,1),\delta _{2}\in (0,1)$ as in \textbf{Condition
4}. Then, from equation (\ref{nevar2bb}),
\begin{eqnarray}
&&\mbox{Var}(S_{2,T})=\sum_{n=m+1}^{\infty }\frac{\mathcal{J}_{n}^{2}}{%
(n!)^{2}}\sigma _{n,K}^{2}(T)  \nonumber
\end{eqnarray}

\begin{eqnarray}
&\leq &k_{1}\left\{ T^{1+\gamma 2d}\left[ \int_{\overline{B}_{T}^{\delta
_{1},\delta _{2}}}+\int_{^{B_{T}^{\delta _{1},\delta _{2}}}}\right] \left( 1-%
\frac{\tau }{T}\right) C^{m+1}(z,\tau )\psi _{T^{\gamma },K}(z)dzd\tau
\right\} ,  \nonumber
\end{eqnarray}%
\noindent for some positive constant $k_{1}>0,$ where we have used $\sum_{n=m}^{\infty}\frac{\mathcal{J}^{2}_{n}}{n!}<\infty,$ and $C^{m+l}\leq
C^{m+1},$ for $l\geq 2,$ since $C(z,\tau )\leq 1,$ under \textbf{Condition 1}.
We then obtain for some positive constants $k_{2},k_{3},k_{4},$
keeping in mind equation (\ref{nevar2bb}), and Remark \ref{remlord}

\begin{eqnarray}
&&\mbox{Var}(S_{2,T})\leq k_{2}\left\{ k_{3}T^{1+\delta _{1}}[T^{\gamma }]^{d(1+\delta _{2})}\right.  \nonumber \\
&&+\left. k_{4} T^{1+\gamma 2d}\sup_{(z,\tau )\in B_{T}^{\delta _{1},\delta
_{2}}}\left\{ C(z,\tau )\right\} \int_{B_{T}^{\delta _{1},\delta
_{2}}}\left( 1-\frac{\tau }{T}\right) C^{m}(z,\tau )\psi _{T^{\gamma
},K}(z)dzd\tau \right\} .  \nonumber
\end{eqnarray}

Hence, as $T\rightarrow \infty ,$
\begin{eqnarray}
&&\frac{\mbox{Var}(S_{2,T})}{\sigma _{m,K}^{2}(T)}\leq k_{5}\left\{ \frac{1}{%
\frac{\sigma _{m,K}^{2}(T)}{T^{1+\delta _{1}}T^{\gamma d(1+\delta _{2})}}}%
+k_{6}\sup_{(z,\tau )\in B_{T}^{\delta _{1},\delta _{2}}}\left\{ C(z,\tau
)\right\} \right.  \nonumber  \label{eqvars2} \\
&&\hspace*{1cm}\left. \hspace*{2cm}\times \frac{\int_{B_{T}^{\delta
_{1},\delta _{2}}}\left( 1-\frac{\tau }{T}\right) C^{m}(z,\tau )\psi
_{T^{\gamma }, K}(z)dzd\tau }{\int_{B_{T}^{\delta _{1},\delta _{2}}\cup
\overline{B}_{T}^{\delta _{1},\delta _{2}}}\left( 1-\frac{\tau }{T}\right)
C^{m}(z,\tau )\psi _{T^{\gamma },K}(z)dzd\tau }\right\} \rightarrow 0,
\nonumber\\
\label{eqc4ii}
\end{eqnarray}%
\noindent where we have applied that from \textbf{Condition 1}, for
all $T>0,$
\begin{equation*}
\frac{\int_{B_{T}^{\delta _{1},\delta _{2}}}\left( 1-\frac{\tau }{T}\right)
C^{m}(z,\tau )\psi _{T^{\gamma },K}(z)dzd\tau }{\int_{B_{T}^{\delta
_{1},\delta _{2}}\cup \overline{B}_{T}^{\delta _{1},\delta _{2}}}\left( 1-%
\frac{\tau }{T}\right) C^{m}(z,\tau )\psi _{T^{\gamma },K}(z)dzd\tau }\leq 1,
\end{equation*}%
\noindent and that from \textbf{Condition 4(i)}, $\sup_{(z,\tau )\in
B_{T}^{\delta _{1},\delta _{2}}}\left\{ C(z,\tau )\right\}\to 0,$ $T\to
\infty,$ as well as from \textbf{Condition 4(ii)} (see also Remark \ref{remlord} and equation (\ref{ineqversus}))
\begin{eqnarray}
&&\frac{\sigma _{m,K}^{2}(T)}{T^{1+\delta _{1}}T^{\gamma d(1+\delta _{2})}}%
\to \infty  \nonumber
\end{eqnarray}
\noindent as $T\to \infty.$ From equations (\ref{eqs1s2}),
(\ref{eqvarfunct}) and (\ref{eqc4ii}),
\begin{equation}
\mbox{Var}\left(Y_{T}-Y_{m,T}\right)\to 0,\quad T\to \infty,
\label{cl2prob}
\end{equation}
\noindent as we wanted to prove.

\begin{remark}
From equation (\ref{cl2prob}), if the limit exists,   $Y_{T}$ and $Y_{m,T}$ have the same limit in probability,    and hence,  in distribution.
\end{remark}

\section{Sojourn functionals for spatiotemporal fields}
\label{exsec}

We consider the first Minkowski functional
\begin{eqnarray}
\hspace*{-3cm} M_{1}(T)&=& \left|\left\{ 0\leq t\leq T,\ x\in T^{\gamma }K;\ Z(x,t)\geq
u\right\}\right|  \notag \\
&=& \int_{0}^{T}\int_{T^{\gamma }K} G_{u}(Z(x,t))dxdt,\quad u\geq 0,\
\gamma\geq 0,  \nonumber
\end{eqnarray}
\noindent where $G_{u}(z)=\mathbb{I}_{z\geq u},$ that is, $G_{u}(Z(x,t))$ is
the indicator function of the set
\begin{equation*}
\left\{0\leq t\leq T, \ x\in T^{\gamma }K; \ Z(x,t)\geq u\right\}.
\end{equation*}
Then,
\begin{eqnarray}
E[M_{1}(T)]&=& (1-\Phi(u))T|K|T^{\gamma d}\nonumber\\
\mathcal{J}_{q}(u)&=& \phi(u)H_{q-1}(u),\quad q\geq 1,\nonumber
\end{eqnarray}
\noindent and from Theorem \ref{th3} with $m=1$ we arrive to the following
result.

\begin{theorem}
\label{MFth1} Under \textbf{Conditions 1--4} with $m=1,$ the random
variable
\begin{eqnarray}
&&\left\{M_{1}(T)-(1-\Phi(u))|K|T^{1+\gamma d}\right\}\times  \notag \\
&&\times \left(\phi (u)\left[2|K|^{2}T^{1+2\gamma d}\int_{0}^{T}\left(1-%
\frac{\tau }{T}\right)\int_{0}^{\mathcal{D}(T^{\gamma
}K)}C(z,\tau)\psi_{T^{\gamma },K}(z)dzd\tau\right]^{1/2}\right)^{-1}  \notag
\end{eqnarray}
\noindent converges to a standard normal distribution as $T\to \infty.$
\end{theorem}

\noindent Note that the case $\gamma =0$ was proved in Theorem 1 of
 \cite{LeonenkoRuiz-Medina2023}. In particular, for
the ball $K=\mathcal{B}(1)\subset \mathbb{R}^{d},$ the random variable
\begin{eqnarray}
&&\left\{M_{1}(T)-(1-\Phi(u))\pi^{d/2}T^{1+\gamma d}\left[\Gamma \left(\frac{d}{2}%
+1\right)\right]^{-1}\right\} \left[\phi (u)\left[(8\pi^{d}\Gamma^{-2}
\left(d/2\right)(1/d))\right.\right.  \notag \\
&&\left.\left.\times T^{1+\gamma d}\int_{0}^{T}\left(1-\frac{\tau }{T}%
\right)\int_{0}^{2T^{\gamma}}z^{d-1}C(z,\tau)I_{1-\left(\frac{z}{2T^{\gamma}}%
\right)^{2}}\left(\frac{d+1}{2},\frac{1}{2}\right)dzd\tau\right]^{1/2}\right]%
^{-1}  \notag
\end{eqnarray}
\noindent is asymptotically distributed as a standard normal random
variable, as $T\to \infty.$

\subsection{Examples}
\label{examples}
Let us analyze \textbf{Condition 4(ii)} with $m=1$ for nonseparable
covariance functions in the Gneiting class (see  \cite{Gneiting02}).
This family of covariance functions is given by

\begin{equation}
\widetilde{C}(\left\Vert x\right\Vert ,\tau )=\frac{\sigma ^{2}}{[\psi (\tau
^{2})]^{d/2}}\varphi \left( \frac{\left\Vert x\right\Vert ^{2}}{\psi (\tau
^{2})}\right) ,\text{ }\sigma ^{2}\geq 0,\text{ }(x,\tau )\in \mathbb{R}%
^{d}\times \mathbb{R},  \label{gneit}
\end{equation}

\noindent in terms of a completely monotone function $\varphi $ and
a positive function $\psi (u),$ $u\geq 0,$ with a completely monotone
derivative.
We will analyze two  special cases of functions $\varphi $ and $%
\psi$  in (\ref{gneit})  in the following two examples.

\subsubsection{Example 1}

 Let us consider the  one--parameter Mittag--Leffler function   $E_{\nu},$ for   $0<\nu \leq 1,$ which
is a completely monotone function (see, e.g.,  \cite{Gorenfloetal2014}; \cite{Haubold}), given by
\begin{equation*}
E_{\nu }(z)=\sum_{k=0}^{\infty}\frac{z^{k}}{
\Gamma (k\beta +1)} ,\quad z\in \mathbb{C},\quad 0<\beta <1.
\end{equation*}

Let function  $\varphi $ in equation  (\ref{gneit})  be defined as
 \begin{equation}
\varphi_{\nu}(z)=E_{\nu }(-z^{\widetilde{\gamma }}),\quad 0<\nu \leq
1,\quad 0< \widetilde{\gamma }<1,  \label{eqmittch}
\end{equation}
\noindent for nonnegative argument $z\geq 0.$ This function  is a complete monotone function, that is,
\begin{equation*}
(-1)^{r}\frac{d^{r}}{dz^{r}}\varphi_{\nu }(z)\geq 0,
\end{equation*}
\noindent for all $r=0,1,2,\dots,$ and $0<\nu \leq 1$ (see  \cite{Barnoff}). From Theorem 4 in  \cite{Simon2014}, for
$z\in \mathbb{R}_{+},$
  that is, for $z\in\mathbb{R}$ and $z\geq 0,$
 and $\nu \in (0,1),$
\begin{eqnarray}
\frac{1}{1+\Gamma (1-\nu )z} &\leq & E_{\nu}(-z)\leq \frac{1}{1+[\Gamma
(1+\nu)]^{-1}z}.  \label{eqimlf}
\end{eqnarray}
The function $\psi (\tau )=(1+a\tau^{\alpha })^{\beta },$ $a>0,$ $0<\alpha \leq 1,$
$0<\beta \leq 1,$ $\tau\geq 0,$ has completely monotone derivatives (see
 \cite{Gneiting02}).

From equations (\ref{b1}) and (\ref{eqimlf}),
\begin{eqnarray}
&& \int_{0}^{T}\left( 1-\frac{\tau }{T}\right) \int_{0}^{2T^{\gamma
}}z^{d-1}C(z,\tau )I_{1-\left( \frac{z}{2T^{\gamma }}\right) ^{2}}\left(
\frac{d+1}{2},\frac{1}{2}\right) dzd\tau  \nonumber \\
&& \geq \int_{0}^{T}\left( 1-\frac{\tau }{T}\right)\int_{0}^{2T^{\gamma }}%
\frac{z^{d-1}\sigma ^{2}}{(a \tau ^{2\alpha }+1)^{\beta d/2}[1+\Gamma
(1-\nu)\left(z^{2\widetilde{\gamma }}/((a \tau ^{2\alpha }+1)^{\beta
\widetilde{\gamma }})\right)] }  \nonumber \\
&& \hspace*{5cm} \times I_{1-\left( \frac{z}{2T^{\gamma }}\right)
^{2}}\left( \frac{d+1}{2},\frac{1}{2}\right) dzd\tau  \nonumber \\
&&=T^{1+\gamma d-\alpha \beta d-2\widetilde{\gamma }\gamma +2 \alpha \beta
\widetilde{\gamma }}\int_{0}^{1}\left( 1-u\right)\int_{0}^{2} \frac{%
x^{d-1}\sigma ^{2}}{(a u ^{2\alpha }+T^{-2\alpha})^{\beta d/2}}  \nonumber \\
&&\hspace*{3cm}\times \frac{I_{1-\left( \frac{x}{2}\right) ^{2}}\left( \frac{%
d+1}{2},\frac{1}{2}\right)}{T^{2 \alpha \beta \widetilde{\gamma }-2%
\widetilde{\gamma }\gamma}+\Gamma (1-\nu)\left[x^{2 \widetilde{\gamma}}/(a u
^{2\alpha }+T^{-2\alpha})^{\beta \widetilde{\gamma }}\right]}dudx, \nonumber\\ \label{chv}
\end{eqnarray}
\noindent where the last equality in equation (\ref{chv}) has been obtained
by applying the change of variables $\tau=Tu$ and $z=T^{\gamma }x.$ From (%
\ref{chv}), \textbf{Condition 4(ii)} holds if
$1+\gamma d-\alpha \beta d-2\gamma \widetilde{\gamma}+2 \alpha \beta
\widetilde{\gamma }>\gamma d \delta_{2}+\delta_{1}$
\noindent for some $\delta_{1}, \delta_{2}\in (0,1).$ In particular, a
sufficient condition for \textbf{Condition 4(ii)} to hold is $\gamma > \alpha \beta  ,$ and $\widetilde{\gamma }< \frac{1}{2(\gamma
-\alpha \beta )}.$
 Under this condition,  one can consider, for instance, $\gamma d-\alpha \beta d>\gamma d \delta_{2},$
and $1-2\gamma \widetilde{\gamma }+2\alpha \beta\widetilde{\gamma }%
>\delta_{1},$ and hence, $0<\delta_{2}<1-\frac{\alpha \beta}{\gamma}<1,$ $%
0<\delta_{1}<1-2\widetilde{\gamma } (\gamma -2\alpha \beta )<1.$
\subsection{Example 2}

 Consider now, in equation (\ref{gneit}),
\begin{eqnarray}
\varphi (z)&=&\frac{1}{(1+\widetilde{c}z^{\widetilde{\gamma }})^{\nu }},\
z>0,\ \widetilde{c}>0,\ 0<\widetilde{\gamma } \leq 1,\ \nu >0  \nonumber \\
\psi (\tau)&=&(1+a\tau^{\alpha })^{\beta },\ a>0,\, 0<\alpha \leq 1,\ 0<\beta \leq
1,\ \tau\geq 0.  \label{fcf}
\end{eqnarray}

It is known  that the function $\varphi (z)$ is completely monotone while, as before,  function $\psi (\tau)$ has completely
monotone derivatives. We restrict our attention to
the ball $\mathcal{B}(1).$ In a similar way to equation (\ref{chv}), one can obtain
\begin{eqnarray}  \label{chv2}
&& \int_{0}^{T}\left( 1-\frac{\tau }{T}\right) \int_{0}^{2T^{\gamma
}}z^{d-1}C(z,\tau )I_{1-\left( \frac{z}{2T^{\gamma }}\right) ^{2}}\left(
\frac{d+1}{2},\frac{1}{2}\right) dzd\tau  \notag \\
&&=T^{1+\gamma d-\alpha \beta d-2\widetilde{\gamma }\gamma \nu+2 \alpha
\beta \widetilde{\gamma }\nu }\int_{0}^{1}\left( 1-u\right)\int_{0}^{2}
\frac{x^{d-1}\sigma ^{2}}{(a u ^{2\alpha }+T^{-2\alpha})^{\beta d/2}}  \notag
\\
&&\hspace*{3cm}\times \frac{I_{1-\left( \frac{x}{2}\right) ^{2}}\left( \frac{%
d+1}{2},\frac{1}{2}\right)}{T^{2 \alpha \beta \widetilde{\gamma }\nu -2%
\widetilde{\gamma }\gamma \nu} +\widetilde{c}\left[x^{2 \widetilde{\gamma}%
\nu }/(a u ^{2\alpha }+T^{-2\alpha})^{\beta \widetilde{\gamma }\nu }\right]}%
dudx.  \notag \\
\end{eqnarray}

\textbf{Condition 4(ii)} is then satisfied if
$ 1+\gamma d-\alpha \beta d-2\gamma \widetilde{\gamma}\nu +2 \alpha \beta
\widetilde{\gamma }\nu >\gamma d \delta_{2}+\delta_{1}, $
 for some $\delta_{1}, \delta_{2}\in (0,1).$ In particular,
a sufficient condition for \textbf{Condition 4(ii)} to hold is
 $\gamma > \alpha \beta ,$ and $\widetilde{\gamma }\nu < \frac{1}{2(\gamma
-\alpha \beta )}. $  Under this condition,  one can consider, for instance,  $\gamma d-\alpha \beta d>\gamma d \delta_{2},$
and $1-2\gamma \widetilde{\gamma }\nu+2\alpha \beta\widetilde{\gamma }
\nu>\delta_{1},$ leading to $0<\delta_{2}<1-\frac{\alpha \beta}{\gamma}<1,$ $
0<\delta_{1}<1-2\nu\widetilde{\gamma } (\gamma -2\alpha \beta )<1.$

%%%%%%%%%%%
\section{Reduction theorem for time varying subordination}
\label{ef}

For each fixed $T>0,$ let $G_{T}\in L_{2}(\mathbb{R},\phi
(z)dz)$ such that
\begin{equation*}
G_{T}(z)=\sum_{q=0}^{\infty }\frac{\mathcal{J}_{q}(T)}{q!}H_{q}(z),\quad
\mathcal{J}_{q}(T)=\int_{\mathbb{R}}G_{T}(z)H_{q}(z)\phi
(z)dz.
\end{equation*}
\noindent For example, one can consider the indicator function with moving threshold,
\begin{equation*}
G_{T}(Z(x,t))=\mathbb{I}_{Z(x,t)\geq u(T)},
\end{equation*}%
\noindent or, equivalently,
\begin{equation*}
G_{T}(z)=\mathbb{I}_{z\geq u(T)}\in L_{2}(\mathbb{R},\phi (z)dz),
\end{equation*}%
\noindent for each fixed $T>0.$ Here, $u:\mathbb{R}_{+}\rightarrow \mathbb{R}
$ be such that $u(T)\rightarrow \infty ,$ as \linebreak $T\rightarrow \infty ,$ and
$G_{T}(z)$ admits the  following orhogonal expansion in terms of Hermite polynomials:
\begin{eqnarray}
\mathcal{J}_{0}(T) &=&\int_{u(T)}^{\infty }\phi (\xi )d\xi =1-\Phi (u(T))
\notag \\
\mathcal{J}_{q}(T) &=&\int_{u(T)}^{\infty }H_{q}(x)\phi
(x)dx=H_{q-1}(u(T))\phi (u(T)),\ q\geq 1.  \label{coefH}
\end{eqnarray}%

  In the subsequent development, we consider the general case of a $T$--varying nonlinear transformation  $G_{T}\in L_{2}(\mathbb{R},\phi
(z)dz).$ We then analyze the asymptotic behavior of the
functional

\begin{eqnarray}
A(T) &=&\int_{0}^{T}\int_{\Lambda (T)K}G_{T}(Z(x,t))dxdt=\sum_{n=0}^{\infty }%
\frac{\mathcal{J}_{n}(T)}{n!}\xi _{n}(T),  \label{f2}
\end{eqnarray}%
where, as before,  $\Lambda (T)=T^{\gamma },\gamma \geq 0,$ and

\begin{equation*}
\xi _{n}(T)=\int_{0}^{T}\int_{T^{\gamma }K}H_{n}(Z(x,t))dxdt,\quad \forall n\in \mathbb{N},\ T>0.
\end{equation*}

The mean and variance  can be computed as follows:
\begin{eqnarray}
E[A(T)]&=&T^{1+\gamma d}|K|\mathcal{J}_{0}(T).  \notag \\
\mbox{Var}(A(T)) &=& \mathbb{E}[A(T)-\mathbb{E}[A(T)]]^{2}=
\sum_{q=m}^{\infty}\frac{\mathcal{J}_{q}^{2}(T)}{(q!)^{2}}%
\sigma^{2}_{q,K\Lambda (T)}(T) , \label{efunvar}
\end{eqnarray}
\noindent where
\begin{eqnarray}
\sigma_{q,\Lambda (T)K}^{2}(T) &=& 2q!T|K|^{2}T^{2d \gamma
}\int_{0}^{T}\left(1-\frac{\tau }{T}\right) \int_{0}^{\mathcal{D}(T^{\gamma
}K)} C^{q}(z,\tau ) \psi_{T^{\gamma }K}(z)dzd\tau.  \label{v2}
\end{eqnarray}
We assume that  function $G_{T}$  has  Hermite rank $m\geq 1,$ for  every $T>0,$ i.e., \textbf{Condition 3} is
satisfied. The derivation of Theorem \ref{th2} below is obtained under the
following additional condition.

\medskip

\noindent \textbf{Condition 5.} For some $m\geq 1,$
\begin{eqnarray}
\overline{\lim}_{T\to \infty}\frac{\mbox{Var}(A(T))}{\frac{\mathcal{J}%
_{m}^{2}(T)}{(m!)^{2}}\sigma^{2}_{m,K\Lambda (T)}(T)}\leq 1.  \label{c5}
\end{eqnarray}

\noindent The technical nature of \textbf{Condition 5}  does not hinder its verification in practice, as follows from the proof of    Theorem \ref{th3b}  below, where this condition is proved to hold for spatiotemporal Minkowski functionals with moving levels.

\begin{theorem}
\label{th2} Under equation (\ref{c5}) in \textbf{Condition 5},

\begin{eqnarray}
Y(T)&=&\frac{\left[A(T)-\mathbb{E}[A(T)]\right]}{|\mathcal{J}%
_{m}(T)|\sigma_{m,T^{\gamma }K}(T)/m!}  \label{funct1} \\
&&  \notag \\
&& \hspace*{-3.75cm}\mbox{and}  \notag \\
&&  \notag \\
&& Y_{m}(T)=\frac{\mbox{sgn}(\mathcal{J}_{m}(T))\int_{0}^{T}\int_{T^{\gamma
}K}H_{m}(Z(x,t))dxdt}{\sigma_{m,T^{\gamma }K}(T)}  \label{funct2}
\end{eqnarray}
\noindent have the same limiting distributions as $T\to \infty$ (if one of
them  exists).
\end{theorem}

\noindent \textbf{Proof.} The proof is straightforward from \textbf{%
Condition 5}. Specifically, from  equation (\ref{efunvar}),
\begin{equation}
\mbox{Var}(A(T))\geq \frac{\mathcal{J}^{2}_{m}(T)\sigma^{2}_{m,\Lambda
(T)K}(T)}{(m!)^{2}}.  \label{in1}
\end{equation}
\noindent Considering $R(T)=Y(T)-Y_{m}(T),$ or, equivalently, $%
Y(T)=R(T)+Y_{m}(T),$ since $\mbox{Var}(Y_{m}(T))=1,$ from equation (\ref{in1}%
), under equation (\ref{c5}),
\begin{equation*}
1=\overline{\lim}_{T\to \infty}\frac{\mbox{Var}(A(T))}{\frac{\mathcal{J}%
_{m}^{2}(T)}{(m!)^{2}}\sigma^{2}_{m,K\Lambda (T)}(T)}=1+ \overline{\lim}%
_{T\to \infty} \mbox{Var}(R(T)).
\end{equation*}
\noindent Thus, $\lim_{T\to \infty}\mbox{Var}(R(T))=0$ (if the limit exists).

\section{First Minkowski functional with moving level}

\label{s5}  We consider the geometric functional
\begin{eqnarray}
M(T)&=& \left|\left\{0\leq t\leq T, \ x\in T^{\gamma }K; \ Z(x,t)\geq
u(T)\right\}\right|  \notag \\
&=&\int_{0}^{T}\int_{T^{\gamma }K} \mathbb{I}_{Z(x,t)\geq u(T)}dxdt  \notag
\\
&=& \left(1-\Phi (u(T))\right)|K| T^{\gamma d}+\sum_{q=1}^{\infty }\frac{%
\mathcal{J}_{q}(T)}{q!} \int_{0}^{T}\int_{T^{\gamma }K} H_{q}(Z(x,t))dxdt,
\label{mf2}
\end{eqnarray}
\noindent where
\begin{equation*}
\mathcal{J}_{q}(T)= \phi(u(T))H_{q-1}(u(T)),
\end{equation*}
\noindent and $u(T)$ is a continuous function such that $u(T)\to \infty,$ as $%
T\to \infty.$

In the next result, the following condition is assumed:

\medskip

\noindent\textbf{Condition 6}. Assume that \textbf{Condition 4} is satisfied, and $u(T)$ is such that $u^{2}(T)=o(log(T)),$ and $u^{2}(T) \sup_{(z,\tau )\in B_{T}^{\beta
_{1},\beta _{2}}}C(z,\tau )\to 0,$ as   $T\to \infty,$ where
\begin{equation}
B_{T}^{\beta _{1},\beta _{2}}=\{(z,\tau );\ \tau \geq T^{\beta _{1}}\ %
\mbox{or}\ z\geq T^{\gamma \beta _{2}}\},  \label{setbeta}
\end{equation}
\noindent for some $\beta_{1}\in (0,\delta_{1}),$ $\beta_{2}\in
(0,\delta_{2}).$

\medskip

\begin{remark}
\label{remordm} Note that, for $$\Lambda (T)=T^{\gamma },\quad  \gamma
\geq 0,\ u^{2}(T)=o(\log(T[\Lambda (T)]^{d}))=o((\gamma d+1)\log (T)),\ T\to \infty,$$
\noindent  if and only if $u^{2}(T)=o(log(T)).$ In particular, for any $\varepsilon_{i}>0,$ $i=1,2,$ $u^{2}(T)=o(\log(T^{\varepsilon_{1}}[\Lambda
(T)]^{d\varepsilon_{2}}))=o((\gamma d \varepsilon_{2}+\varepsilon_{1})\log
(T))=o(\log (T)).$
\end{remark}

\bigskip

\begin{remark}
Note that the set $B_{T}^{\delta_{1},\delta_{2}}$ introduced in equation  (\ref{cz}) is not included in the set family  $\left\{B_{T}^{\beta _{1},\beta_{2}},\ \beta_{1}\in (0,\delta_{1}),\ \beta_{2}\in (0,\delta_{2})\right\}.$  This set family is considered in the proof of Theorem  \ref{th3b} (see equations (\ref{cal2})--(\ref{s1infty})), to apply \textbf{Condition 6} incorporating \textbf{Condition 4(ii)} (see also Remark \ref{remordm}).
\end{remark}

\bigskip

\begin{theorem}
\label{th3b} Under \textbf{Conditions 1--3}, and \textbf{Condition 6}, as  $T\rightarrow \infty ,$ the random variables
\begin{equation*}
\frac{M(T)-T^{1+\gamma d}|K|(1-\Phi (u(T)))}{\phi (u(T))\left[
2T|K|^{2}T^{2\gamma d}\int_{0}^{T}\left( 1-\frac{\tau }{T}\right) \int_{0}^{%
\mathcal{D}(T^{\gamma }K)}C(z,\tau )\psi _{\rho _{T^{\gamma }K}}(z,\tau
)dzd\tau \right] ^{1/2}}
\end{equation*}%
\noindent and
\begin{equation*}
\frac{\int_{0}^{T}\int_{T^{\gamma }K}Z(x,t)dxdt}{\left[ 2T|K|^{2}T^{2d\gamma
}\int_{0}^{T}\left( 1-\frac{\tau }{T}\right) \int_{0}^{\mathcal{D}(T^{\gamma
}K)}C(z,\tau )\psi _{\rho _{T^{\gamma }K}}(z,\tau )dzd\tau \right] ^{1/2}}
\end{equation*}%
\noindent have the same asymptotic distribution, that is, a standard normal
distribution, where $M(T)$ has been introduced in (\ref{mf2}).
\end{theorem}

\noindent \textbf{Proof.} It is known that for bivariate normal
density (see, e.g., equation 10.8.3 in  \cite{CramerLeadbetter})
\begin{eqnarray}
&& \phi(x,y,\rho)= \frac{1}{2\pi\sqrt{1-\rho^{2}}}\exp\left(-\frac{1}{%
2(1-\rho^{2})}(x^{2}+y^{2}-2\rho x y)\right),  \label{eq1}
\end{eqnarray}
\noindent the following identity holds:
\begin{eqnarray}
&&\int_{u(T)}^{\infty}\int_{u(T)}^{\infty}\phi(x,y,\rho)dxdy  \notag \\
&&=\left(\int_{u(T)}^{\infty}\phi(y)dy\right)^{2}+\frac{1}{2\pi}
\int_{0}^{\rho}\exp\left(-\frac{u^{2}(T)}{1+v}\right)\frac{dv}{\sqrt{1-v^{2}}%
}.  \label{ednd}
\end{eqnarray}

\noindent In our case, in equation (\ref{ednd}), $\rho=\widetilde{C}(\|x-y\|,|t-s|),$   we then obtain
\begin{eqnarray}  \label{eqf1}
&&E[M(T)]=(1-\Phi(u(T)))|K|T^{1+\gamma d}\nonumber\\
&& E[M^{2}(T)]=\int_{[0,T]\times [0,T]}\int_{T^{\gamma }K\times T^{\gamma }K}
\mathbb{E}\left[\mathbb{I}_{Z(x,t)>u(T)}\mathbb{I}_{Z(y,s)>u(T)}\right]
dxdydtds  \nonumber \\
&&=\int_{[0,T]\times [0,T]}\int_{T^{\gamma }K\times T^{\gamma
}K}\int_{u(T)}^{\infty}\int_{u(T)}^{\infty}\phi(u,w,\rho)dudwdxdydtds  \notag
\\
&&= T^{2+2d \gamma }|K|^{2}\left[ \int_{u(T)}^{\infty}\phi(u)du\right]^{2}
\nonumber  \\
&&+ \frac{1}{2\pi}\int_{0}^{T}\int_{0}^{T}\int_{T^{\gamma }K\times T^{\gamma
}K} \int_{0}^{\widetilde{C}(\|x-y\|,|t-s|)}\exp\left(-\frac{u^{2}(T)}{1+v}%
\right)\frac{dv}{\sqrt{1-v^{2}}}dxdydtds  \nonumber  \\
&&= T^{2+2d \gamma }|K|^{2}[1-\Phi (u(T))]^{2}  \notag \\
&&+\frac{1}{2\pi} \int_{0}^{T}\int_{0}^{T}\int_{T^{\gamma }K\times T^{\gamma
}K}\int_{0}^{\widetilde{C}(\|x-y\|,|t-s|)}\exp\left(-\frac{u^{2}(T)}{1+v}%
\right)\frac{dv}{\sqrt{1-v^{2}}}dxdydtds.\nonumber  \\
\end{eqnarray}
\noindent Thus, from equation (\ref{eqf1}), operating in a similar way to equation (\ref{nevar2bb}) in the integrals over time and space,  we have
\begin{eqnarray}
\mbox{Var}(M(T))&=& E\left[M^{2}(T)\right]-\left[E[M(T)]\right]^{2}
\nonumber  \\
&=&
2T|K|^{2}T^{2\gamma d}\int_{0}^{T}\left(1-\frac{\tau}{T}%
\right)\int_{0}^{\mathcal{D}(T^{\gamma }K)}z^{d-1}\widetilde{\psi}%
_{\rho_{T^{\gamma }K}}(z)  \notag \\
& & \hspace*{2cm} \times\frac{1}{2\pi}\int_{0}^{C(z,\tau)}\exp\left(-\frac{u^{2}(T)}{1+v}%
\right)\frac{dv}{\sqrt{1-v^{2}}}dzd\tau,   \label{vmf}
\end{eqnarray}
\noindent where $\widetilde{\psi}_{\rho_{T^{\gamma
}K}}(z)=(z^{-(d-1)})\psi_{\rho_{T^{\gamma }K}}(z),$ with $%
\psi_{\rho_{T^{\gamma }K}}(z)$ being, as before, the probability density of
the random variable $\rho_{K}=\|P_{1}-P_{2}\|,$ where $P_{1}$ and $P_{2}$
are two independent random points with uniform distribution in $T^{\gamma }K$
(see also Remark \ref{remlord}). The proof is based on verifying that
\textbf{Condition 5} holds under the assumptions made, and hence, applying
Theorem \ref{th2} for $m=1.$ Specifically, we prove that

\begin{equation*}
\overline{\lim}_{T\to \infty}\frac{\mbox{Var}(M(T))}{\mathcal{J}%
_{1}^{2}(T)\sigma^{2}_{1,KT^{\gamma }}(T)}\leq 1.
\end{equation*}
Note that $\mathcal{J}_{1}(T)=\phi(u(T))=\frac{1}{\sqrt{2\pi}}%
\exp(-(u^{2}(T))/2)>0,$ and
\begin{eqnarray}
\frac{\exp\left(-\frac{u^{2}(T)}{1+v}\right)}{\phi^{2}(u(T))}&=& 2\pi
\exp\left(u^{2}(T)\left(\frac{v}{1+v}\right)\right).  \label{fv}
\end{eqnarray}
Hence, from equations (\ref{vmf}) and (\ref{fv}),
\begin{eqnarray}
&& \frac{\mbox{Var}(M(T))}{[\phi(u(T))]^{2}}= 2 T|K|^{2}T^{\gamma
2d}\int_{0}^{T}\left(1-\frac{\tau }{T}\right) \int_{0}^{\mathcal{D}%
(T^{\gamma }K)}z^{d-1}\widetilde{\psi}_{\rho_{T^{\gamma }K}}(z)  \notag \\
&& \hspace*{2cm}\times \int_{0}^{C(z,\tau)}\exp\left(u^{2}(T)\left(\frac{v}{%
1+v}\right)\right)\frac{dv}{\sqrt{1-v^{2}}}dzd\tau  \notag \\
&&=2 T|K|^{2}T^{2\gamma d}\left[\int_{\overline{B}_{T}^{\beta_{1},\beta_{2}
}}+\int_{B_{T}^{\beta_{1},\beta_{2}}}\right]\left(1-\frac{\tau }{T}\right)
z^{d-1}\widetilde{\psi}_{\rho_{T^{\gamma }K}}(z)  \notag \\
&& \hspace*{2cm}\times \int_{0}^{C(z,\tau)}\exp\left(u^{2}(T)\left(\frac{v}{%
1+v}\right)\right)\frac{dv}{\sqrt{1-v^{2}}}dzd\tau  \notag \\
&&= S_{1}(T)+S_{2}(T),  \label{cal}
\end{eqnarray}
\noindent where the set $B_{T}^{\beta _{1},\beta _{2}}$ has been introduced
in equation (\ref{setbeta}), and the set $\overline{B}_{T}^{\beta_{1},%
\beta_{2}}$ is defined as
\begin{equation*}
\overline{B}_{T}^{\beta_{1},\beta_{2}}=\left\{ (z,\tau);\ 0\leq \tau \leq
T^{\beta_{1}},\ 0\leq z\leq T^{\gamma \beta_{2}}\right\},
\end{equation*}
\noindent with $\beta_{1} \in (0,\delta_{1}),$ $\beta_{2} \in
(0,\delta_{2}). $

From \textbf{Condition 6}, applying Remark \ref{remordm}, $u^{2}(T)=o\left(\log((T)^{\varepsilon_{1}}([\Lambda
(T)]^{d})^{\varepsilon_{2}})\right),$ for any $\varepsilon_{1},\varepsilon_{2}>0.$ Keeping in mind that $\frac{v}{1+v}\leq 1,$ we obtain

\begin{eqnarray}
&& S_{1}(T) \leq k_{8} \ T^{1+2\gamma d}\int_{\overline{B}%
_{T}^{\beta_{1},\beta_{2} }}\left(1-\frac{\tau }{T}\right)z^{d-1}\widetilde{%
\psi}_{\rho_{T^{\gamma }K}}(z)  \notag \\
&&\hspace*{2cm}\times \exp\left( u^{2}(T)\right)\int_{0}^{C(z,\tau)}\frac{dv%
}{\sqrt{1-v^{2}}}dzd\tau  \notag \\
& &= k_{8}\ T^{1+2\gamma d}\exp\left( u^{2}(T)\right)\int_{\overline{B}%
_{T}^{\beta_{1},\beta_{2} }}\left(1-\frac{\tau }{T}\right) z^{d-1}\widetilde{%
\psi}_{\rho_{T^{\gamma }K}}(z)  \notag \\
&&\hspace*{2cm}\times \mbox{arc sin}C(z,\tau)dzd\tau  \notag \\
& &\leq \widetilde{k_{8}} \ T^{1+\varepsilon_{1}}T^{2d \gamma +d \gamma
\varepsilon_{2}}\int_{\overline{B}_{T}^{\beta_{1},\beta_{2} }}\left(1-\frac{%
\tau }{T}\right)z^{d-1}\widetilde{\psi}_{\rho_{T^{\gamma }K}}(z)dzd\tau
\notag \\
& &\leq k_{9} T^{1+\beta_{1} +\varepsilon_{1}}[\Lambda (T)]^{d+d\beta_{2}
+d\varepsilon_{2}},  \label{cal2}
\end{eqnarray}
\noindent for some positive constants $k_{8},$ and  $ k_{9}.$  In particular, for $%
\varepsilon_{i} =\delta_{i}-\beta_{i},$ $i=1,2,$ under
\textbf{Condition 4(ii)}, as $T\to \infty,$
\begin{eqnarray}
&& \frac{S_{1}(T)}{\sigma_{1,K}^{2}(T)}\leq \frac{k_{9}}{\frac{1}{%
T^{1+\delta_{1} }T^{\gamma d(1+\delta_{2})}}\sigma_{1,K}^{2}(T)}\to 0.
\label{s1infty}
\end{eqnarray}

Let us now consider $S_{2}(T),$ under  \textbf{Condition 6},
\begin{eqnarray}
&&\sup_{(z,\tau )\in B_{T}^{\beta _{1},\beta _{2}}}C(z,\tau )\rightarrow
0,\quad T\rightarrow \infty  \notag \\
&&B_{T}^{\beta_{1},\beta _{2}}=\{(z,\tau );\ \tau \geq T^{\beta _{1}}\ %
\mbox{or}\ z\geq T^{\gamma \beta_{2}}\},  \nonumber
\end{eqnarray}

\noindent and,  since  $\frac{1}{1+v}\leq 1$ and $\frac{1}{\sqrt{1-v^{2}}}
\to 1,$  $0\leq v\leq C(z,\tau),$ $C(z,\tau)\to 0,$ $(z,\tau)\in B_{T}^{\beta _{1},\beta _{2}},$ as $T\to \infty,$
  we obtain  that $S_{2}(T)$ satisfies
\begin{eqnarray}
&& S_{2}(T)= 2T^{1+2\gamma d}|K|^{2}\int_{B_{T}^{\beta_{1},\beta_{2}
}}\left(1-\frac{\tau}{T}\right) \psi_{\rho_{T^{\gamma }K}}(z)  \nonumber \\
& &\hspace*{3.5cm} \times \int_{0}^{C(z,\tau)}\exp\left(u^{2}(T)\frac{v}{1+v}%
\right)\frac{dv}{\sqrt{1-v^{2}}}dzd\tau  \nonumber \\
& & \leq 2T^{1+2\gamma d}|K|^{2}\int_{B_{T}^{\beta_{1},\beta_{2}}}\left(1-%
\frac{\tau}{T}\right) \psi_{\rho_{T^{\gamma }K}}(z)  \nonumber \\
& &\hspace*{3.5cm} \times \int_{0}^{C(z,\tau)}\exp\left(u^{2}(T)v\right)%
\frac{dv}{\sqrt{1-v^{2}}}dzd\tau  \nonumber \end{eqnarray}

 \begin{eqnarray}
& &\leq k_{10} 2T^{1+2\gamma d}|K|^{2}\exp\left(u^{2}(T) \sup_{(z,\tau )\in
B_{T}^{\beta _{1},\beta _{2}}}C(z,\tau ) \right)  \nonumber\\
& & \hspace*{3.5cm}\times \int_{B_{T}^{\beta_{1},\beta_{2} }}\left(1-\frac{%
\tau}{T}\right) \psi_{\rho_{T^{\gamma }K}}(z) C(z,\tau)dzd\tau ,  \nonumber
\end{eqnarray}
\noindent for some positive constant $k_{10}.$ Hence, for $T$ sufficiently
large,
\begin{eqnarray}
&&\frac{S_{2}(T)}{\sigma_{1,K}^{2}(T)}\leq k_{11}\exp\left(u^{2}(T)
\sup_{(z,\tau )\in B_{T}^{\beta _{1},\beta _{2}}}C(z,\tau ) \right)  \notag
\\
& & \hspace*{1.5cm} \times \frac{\int_{B_{T}^{\beta_{1},\beta_{2}}}\left(1-%
\frac{\tau}{T}\right) \psi_{\rho_{T^{\gamma }K}}(z)C(z,\tau)dzd\tau }{\int_{
\overline{B}_{T}^{\beta_{1},\beta_{2} } \cup B_{T}^{\beta_{1},\beta_{2}
}}\left(1-\frac{\tau}{T}\right) \psi_{\rho_{T^{\gamma }K}}(z)C(z,\tau)dzd\tau%
}\leq 1,  \label{ee}
\end{eqnarray}
\noindent since under \textbf{Condition 6},
\begin{equation*}
\exp\left(u^{2}(T) \sup_{(z,\tau )\in B_{T}^{\beta _{1},\beta _{2}}}C(z,\tau
) \right)\to 1,\quad T\to \infty,
\end{equation*}
\noindent and $k_{11}\leq 1,$ for sufficiently large $T.$ Finally, from (\ref{s1infty})--(\ref{ee}),
\begin{eqnarray}
&& \vspace*{-3cm} \overline{\lim}_{T\to \infty}\frac{\mbox{Var}(M(T))}{%
\phi^{2}(u(T))\sigma_{1, K}^{2}(T)}  \notag \\
&&= \overline{\lim}_{T\to \infty}\left(\frac{S_{1}(T)}{\sigma_{1, K}^{2}(T)}+%
\frac{S_{2}(T)}{\sigma_{1, K}^{2}(T)}\right)\leq 1.  \label{elimfv}
\end{eqnarray}
\noindent Thus, the desired result follows from Theorem \ref{th2}, and the asymptotic normality of \begin{equation*}
\frac{\int_{0}^{T}\int_{T^{\gamma }K}Z(x,t)dxdt}{\left[ 2T|K|^{2}T^{2d\gamma
}\int_{0}^{T}\left( 1-\frac{\tau }{T}\right) \int_{0}^{\mathcal{D}(T^{\gamma
}K)}C(z,\tau )\psi _{\rho _{T^{\gamma }K}}(z,\tau )dzd\tau \right] ^{1/2}}.
\end{equation*}
\section{Sojourn functionals for a class of spherical random fields}
\label{sec7}
This section derives a central limit result for sojourn functionals
subordinated to   STGRFs homogeneous and isotropic in space, and stationary in time,
 restricted to the unit sphere.

 Let $\mathbb{S}_{d-1}(1)=\{x\in \mathbb{R}^{d};\ \Vert x\Vert =1\}$ be
the unit sphere embedded into $\mathbb{R}^{d},$ for some $d\geq 2,$ and
denote by $d\nu _{d-1}(x)$ the normalized Riemannian measure on $\mathbb{S}%
_{d-1}(1).$ Denote also by $\theta =\arccos \left( \left\langle x,x^{\prime
}\right\rangle \right) $ the angle between two points $x,x^{\prime }\in
\mathbb{S}_{d-1}(1).$  For every $x,x^{\prime}\in \mathbb{S}_{d-1}(1),$   $\Vert
x-x^{\prime }\Vert =2\sin \left( \frac{\theta }{2}\right),$ with $\|\cdot\|$  being the Euclidean distance. Let us denote
$S_{lm}^{(d)}(u),$ $u\in \mathbb{S}_{d-1}(1),$ $m=1,2,\dots ,h(l,d),$ $l\in
\mathbb{N}_{0},$ the real spherical harmonics on $\mathbb{S}_{d-1}(1)$ (see
 \cite{Leonenko99};  \cite{MullerC}, and the references
therein), with $h(l,d)$\linebreak $=(2l+d-2)\frac{(l+d-3)!}{(d-2)!l!}$ denoting the
dimension of the eigenspace of the Laplace Beltrami operator  generated by

$$\left\{S_{lm}^{(d)},
\quad  m=1,2,\dots ,h(l,d)\right\},\quad  l\in \mathbb{N}_{0}.$$

Let $\left\{ Z(x,t),\ x\in \mathbb{R}^{d},\ t\in \mathbb{R}\right\}$
be a zero--mean,   mean--square continuous,
 STGRF,  homogeneous and  isotropic in space, and stationary in time, with covariance function satisfying
\begin{eqnarray}
&&\widetilde{C} (\|x\|,\tau )  =2^{\frac{d-2}{2}+1}\Gamma \left(\frac{d}{2}\right)\int_{0}^{\infty}\cos%
\left(\mu \tau\right) \int_{0}^{\infty}\frac{J_{\frac{d-2}{2}}\left(\lambda
\|x\|\right)}{\left(\lambda \|x\|\right)^{\frac{d-2}{2}}}\mathcal{G}(d\lambda, d\mu ),\nonumber\\
\label{bimeasure}
\end{eqnarray}
\noindent where $\mathcal{G}$ is defined from the spectral measure $F$ of $Z$ arising  in Bochner Theorem (see, e.g., \cite{Ivanov89}; \cite{Schoenberg38})  satisfying

\begin{eqnarray}
\int_{\mathbb{R}}\int_{\mathbb{R}^{d}}F(d\widetilde{\omega },d\widetilde{%
\mu })&=&\int_{0}^{\infty}\int_{0}^{\infty} \mathcal{G}(d\lambda , d\mu)<\infty,  \notag
\\
 \mathcal{G}(\lambda ,\mu )&=&\int_{\|\widetilde{\omega }\|< \lambda }\int_{|%
\widetilde{\mu }|<\mu}F(d\widetilde{\omega }, d\widetilde{\mu }). \label{sm}
\end{eqnarray}

Let us consider  $T_{R}(x,t)= Z(x,t),$ for every $x\in \mathbb{S}_{d-1}(1),$ and  $t\in \mathbb{R},$ defining  the restriction of  $Z(x,t)$ to the unit sphere $\mathbb{S}_{d-1}(1).$
The following identities will be applied in the characterization of the second--order pure point  spectral  properties of $T_{R}(x,t)$ from the spectral representation
(\ref{bimeasure}) of the covariance function $\widetilde{C}$ of $Z.$

 In the following, we denote
  $\mathbb{S}_{d-1}(u)=\left\{x\in \mathbb{R}^{d}; \ \|x\|=u\right\}.$
 Let us first consider the characteristic function of  $\mathbb{S}_{d-1}(u)$
\begin{eqnarray}
&& \frac{1}{|\mathbb{S}_{d-1}(u)|}\int_{\mathbb{S}_{d-1}(u)}\exp\left(i\left%
\langle \lambda ,x\right\rangle\right)d\nu_{d-1} (\lambda )=Y_{d}(ux),
\label{rssph}
\end{eqnarray}
\noindent where $Y_{d}$ denotes the spherical Bessel function, and $d\nu_{d-1} $
is the normalized Riemannian measure on $\mathbb{S}_{d-1}(1).$
Applying equation (1.2.13a) in \cite{Ivanov89}, the following relationship
holds between the spherical Bessel function and the Bessel function of the
first kind
\begin{equation}
Y_{d}(z)= 2^{(d-2)/2}\Gamma \left(\frac{d}{2}\right) J_{\frac{d-2}{2}%
}(z)z^{(2-d)/2},\quad z\geq 0.  \label{spbf}
\end{equation}
\noindent Here,
\begin{equation*}
J_{\nu }(z)=\sum_{m=0}^{\infty}(-1)^{m}\left(\frac{z}{2}\right)^{2m+\nu}%
\left[m!\Gamma (m+\nu +1)\right]^{-1}
\end{equation*}
\noindent is the Bessel function of the first kind of order $\nu>1/2.$

Thus, from equations (\ref{rssph}) and (\ref{spbf}),
\begin{eqnarray}
&& \hspace*{-2cm}\frac{1}{|\mathbb{S}_{d-1}(u)|}\int_{\mathbb{S}_{d-1}(u)}\exp\left(i\left%
\langle \lambda, x-x^{\prime }\right\rangle\right)d\nu _{d-1} (\lambda )  \notag \\
&&=2^{(d-2)/2}\Gamma \left(\frac{d}{2}\right)\frac{J_{\frac{d-2}{2}%
}(u\|x-x^{\prime }\|)}{\left(u\|x-x^{\prime }\|\right)^{\frac{d-2}{2}}}.
\label{bffk}
\end{eqnarray}

Applying now addition theorem of spherical Bessel function
\begin{equation}
Y_{d}(\lambda \rho )=c_{1}^{2}(d)\sum_{l=0}^{\infty
}\sum_{m=1}^{h(l,d)}S_{lm}^{(d)}(u)S_{lm}^{(d)}(v)\frac{J_{l+\frac{d-2}{2}%
}(\lambda r_{1})}{(\lambda r_{1})^{\frac{d-2}{2}}}\frac{J_{l+\frac{d-2}{2}%
}(\lambda r_{2})}{(\lambda r_{2})^{\frac{d-2}{2}}},  \label{eqaddf}
\end{equation}%

\noindent where $c_{1}^{2}(d)=2^{d-1} \Gamma \left(\frac{d}{2}%
\right) \pi^{d/2},$ and
\begin{eqnarray}  \label{eqsphcoor}
x&=& (r_{1},u),\ r_{1}\geq 0,\ u=x/\|x\|\in \mathbb{S}_{d-1}(1)  \notag \\
y&=&(r_{2},v),\ r_{2}\geq 0,\ v=y/\|y\|\in \mathbb{S}_{d-1}(1)  \notag \\
\rho &=& \|x-y\|=\sqrt{r_{1}^{2}+r_{2}^{2}-2r_{1}r_{2}\cos (\gamma )},\quad
\cos (\gamma )=\frac{\left\langle x,y \right\rangle}{\|x\|\|y\|},\ \lambda
\geq 0.  \notag \\
\end{eqnarray}

From equations (\ref{rssph})--(\ref{eqsphcoor}),

\begin{eqnarray}
&&\widetilde{C} (\|x-y\|,\tau)  \notag \\
&&=2[c_{1}(d)]^{2} \sum_{l=0}^{\infty}\left[\int_{0}^{\infty}\int_{0}^{%
\infty}\cos(\mu \tau)\frac{J_{l+\frac{d-2}{2}}(\lambda r_{1})}{(\lambda
r_{1})^{\frac{d-2}{2}}}\frac{J_{l+\frac{d-2}{2}}(\lambda r_{2})}{(\lambda
r_{2})^{\frac{d-2}{2}}}\mathcal{G}(d\lambda ,d\mu)\right]  \notag \\
&& \hspace*{5cm} \times\sum_{m=1}^{h(l,d)} S_{lm}^{(d)}(u)S_{lm}^{(d)}(v),
\label{eqexorth}
\end{eqnarray}
\noindent where $(x,y),$ $(u,v),$ and $(r_{1},r_{2})$ are defined as in
equation (\ref{eqsphcoor}). For $r_{1}=r_{2}=1,$ that is, in the
case of  considering the covariance function $\widetilde{C}_{R}$ of the restricted random field $\left\{ T_{R}(x,t),\ x\in \mathbb{S}_{d-1}(1), \ t\in \mathbb{R}\right\},$  we obtain
\begin{eqnarray}  \label{eqexorth}
&&\widetilde{C}_{R} (\|x-y\|,\tau )  \notag \\
&&=2[c_{1}(d)]^{2} \sum_{l=0}^{\infty}\left[\int_{0}^{\infty}\int_{0}^{%
\infty}\cos(\mu \tau)\left[\frac{J_{l+\frac{d-2}{2}}(\lambda )}{\lambda^{%
\frac{d-2}{2}}}\right]^{2}\mathcal{G}(d\lambda ,d\mu)\right]  \notag \\
&&\hspace*{3.5cm}\times
\sum_{m=1}^{h(l,d)}S_{lm}^{(d)}(u)S_{lm}^{(d)}(v).
\end{eqnarray}

\noindent Random field  $\left\{ T_{R}(x,t),\ x\in \mathbb{S}_{d-1}(1), \ t\in \mathbb{R}\right\}$
 then has $\tau$--varying angular power spectrum $\left\{ A_{l}(\tau
),\ \tau \geq 0,\ l\in \mathbb{N}_{0}\right\}$ given by

\begin{eqnarray}  \label{eqas}
A_{l}(\tau )=2\Gamma \left( \frac{d}{2}\right)\pi^{d/2}\int_{0}^{\infty}%
\int_{0}^{\infty}\cos(\mu \tau)\left[\frac{J_{l+\frac{d-2}{2}}(\lambda )}{%
\lambda^{\frac{d-2}{2}}}\right]^{2}\mathcal{G}(d\lambda ,d\mu),\ \tau\geq 0, \ l\in
\mathbb{N}_{0}.  \notag \\
\end{eqnarray}

Equation (\ref{eqas}) allows  the interpretation of the elements of the
$\tau$--varying angular spectrum $\left\{ A_{l}(\tau ),\ \tau\geq 0,\ l\in \mathbb{N}%
_{0}\right\}$ as the inverse Fourier transforms of the temporal spectral
measures $f(d\mu )=2\Gamma \left(\frac{d}{2}\right)\pi^{d/2}\int_{0}^{\infty}%
\left[\frac{J_{l+\frac{d-2}{2}}(\lambda )}{\lambda^{\frac{d-2}{2}}}\right]%
^{2}\mathcal{G}(d\lambda ,d\mu).$ Equivalently,

\begin{eqnarray}
A_{l}(\tau )&=&\int_{0}^{\infty }\cos (\mu \tau )f(d\mu )\nonumber\\ &=&2\Gamma \left( \frac{%
d}{2}\right) \pi ^{d/2}\int_{0}^{\infty }\cos (\mu \tau )\int_{0}^{\infty }%
\left[ \frac{J_{l+\frac{d-2}{2}}(\lambda )}{\lambda ^{\frac{d-2}{2}}}\right]
^{2}\mathcal{G}(d\lambda ,d\mu ).
\label{timevaryas}
\end{eqnarray}

Note also that, from  (\ref{bimeasure}), applying trigonometric identity $\Vert
x-x^{\prime }\Vert =2\sin \left( \frac{\theta }{2}\right),$ for every $x,x^{\prime}\in \mathbb{S}_{d-1}(1),$
\begin{eqnarray}
&&\widetilde{C}_{R} (\|x-x^{\prime}\|,\tau )  =
2^{\frac{d-2}{2}+1}\Gamma \left(\frac{d}{2}\right)\int_{0}^{\infty}\cos%
\left(\mu \tau\right) \int_{0}^{\infty}\frac{J_{\frac{d-2}{2}}\left(\lambda
\|x\|\right)}{\left(\lambda \|x\|\right)^{\frac{d-2}{2}}}\mathcal{G}(d\lambda, d\mu )
\nonumber\\
&& 2^{\frac{d-2}{2}+1}\Gamma \left(\frac{d}{2}\right)\int_{0}^{\infty}\cos%
\left(\mu \tau\right) \int_{0}^{\infty}\frac{J_{\frac{d-2}{2}}\left(\lambda
2\sin \left( \frac{\theta }{2}\right)\right)}{\left(\lambda 2\sin \left( \frac{\theta }{2}\right)\right)^{\frac{d-2}{2}}}\mathcal{G}(d\lambda, d\mu )\nonumber\\ &&=C_{R}(\cos (\theta ),\tau ).
\label{bimeasureR}
\end{eqnarray}
\noindent  Here,  as before, $\theta $ denotes the angle between vectors $x$ and $x^{\prime}$ in $\mathbb{S}_{d-1}(1).$
Thus, $T_{R}$  has covariance function (\ref{bimeasureR}) (see, e.g., \cite{Leonenkoetal17a}).

 \medskip

 \noindent \emph{Condition K.}  Random field  $\left\{T_{R}(x,t),\
x\in \mathbb{S}_{d-1}(1),\ t\in \mathbb{R}\right\}$ is defined on the sphere as the restriction of a zero--mean
 STGRF
with covariance function  (\ref{bimeasure}).

\medskip

Under \emph{Condition K} (see equations (\ref{eqexorth})--(\ref{timevaryas})),   random
field $T_{R}$ admits the following orthogonal expansion, in the mean--square sense,   for every fixed $t\in \mathbb{R},$  and $x\in \mathbb{S}_{d-1}(1),$ $$T_{R}(x,t)=\sum_{l=0}^{\infty}\sum_{m=1}^{h(l,d)}a_{lm}(t)S_{lm}^{(d)}(x),$$
\noindent   where $a_{lm}(t),$ $m=1,\dots, h(l,d),$ $l\in \mathbb{N}_{0},$ are independent zero--mean Gaussian stochastic processes such that $E[a_{lm}(t)]=0,$
$E[a_{lm}(t)a_{l^{\prime }m^{\prime}}(t^{\prime })]=\delta_{ll^{\prime }}\delta_{mm^{\prime }}A_{l}(|t-t^{\prime }|),$ with $A_{l}(|t-t^{\prime }|) =A_{l}(\tau )$ satisfying  (\ref{timevaryas}), and
$\sum_{l=0}^{\infty }h(l,d)A_{l}(\tau )<\infty,$ for every $\tau \in \mathbb{R}_{+}.$

 Let now consider the first Minkowski
functional subordinated to $T_{R}(x,t),$ given by

\begin{eqnarray}
&&\hspace*{-2cm}N_{T}=\int_{0}^{T}\int_{\mathbb{S}_{d-1}(1)}\mathbb{I}_{T_{R}(x,t)\geq u}d\nu
_{d-1}(x)dt  \notag \\
&=&\left\vert \left\{ 0\leq t\leq T;\ T_{R}(x,t)\geq u,\ x\in \mathbb{S}_{d-1}(1)
\right\} \right\vert   \notag \\
&=&E[N_{T}]+\sum_{n\geq 1}\frac{\mathcal{J}_{n}}{n!}\eta_{n}(T),
\label{fmftr}
\end{eqnarray}%
\noindent where $E[N_{T}]=(1-\Phi (u))\frac{2\pi ^{d/2}}{\Gamma \left( \frac{%
d}{2}\right) },$ and $\mathcal{J}_{n}(u)=\phi (u)H_{n-1}(u),$ $n\geq 1,$ and
\begin{equation*}
\eta _{n}(T)=\int_{0}^{T}\int_{\mathbb{S}_{d-1}(1)}H_{n}(T_{R}(x,t))d\nu
_{d-1}(x)dt.
\end{equation*}

\noindent Thus, $E[\eta _{n}(T)]=0,$ $E[\eta _{n}(T)\eta _{l}(T)]=0,$
$n\neq l,$ and
\begin{eqnarray}
&&\sigma _{n}^{2}(T)=E[\eta _{n}^{2}(T)]  \notag \\
&=&2n!\int_{0}^{T}\int_{0}^{T}\int_{\mathbb{S}_{d-1}(1)\times \mathbb{S}_{d-1}(1)}%
\widetilde{C}^{n}\left( \Vert x-x^{\prime }\Vert ,|t-t^{\prime }|\right)
d\nu (x)d\nu (x^{\prime })dtdt^{\prime }  \notag \\
&=&2n!T|\mathbb{S}_{d-1}(1)|^{2}\int_{0}^{T}\left( 1-\frac{\tau }{T}\right)
E\left( \widetilde{C}^{n}(\Vert W_{1}-W_{2}\Vert ,\tau )\right) d\tau
\notag \\
&=&2n!T|\mathbb{S}_{d-1}(1)|^{2}\frac{1}{\sqrt{\pi }}\Gamma \left( \frac{d}{2}%
\right) \Gamma ^{-1}\left( \frac{d-1}{2}\right)\nonumber\\
&&\hspace*{1cm} \times
\int_{0}^{T}\int_{0}^{2}\left( 1-\frac{\tau }{T}\right) z^{d-2}\left( 1-%
\frac{z^{2}}{4}\right) ^{\frac{d-3}{2}}C^{n}(z,\tau )dzd\tau\nonumber\\
&=&
 2n!T4\pi^{d-1/2}\Gamma \left( \frac{d}{2}%
\right) \left[\Gamma \left( \frac{d-1}{2}\right)\right]^{-1}\nonumber\\
&&\hspace*{1cm} \times
\int_{0}^{T}\int_{0}^{2}\left( 1-\frac{\tau }{T}\right) z^{d-2}\left( 1-%
\frac{z^{2}}{4}\right) ^{\frac{d-3}{2}}C^{n}(z,\tau )dzd\tau ,
\label{eqchaoscomp}
\end{eqnarray}
\noindent where  $C^{n}$ denotes, as before, the $n$th power of the covariance function $C.$
Here, $W_{1}$ and $W_{2}$ are two independent uniformly
distributed random vectors on $\mathbb{S}_{d-1}(1)$ with probability density of
their Euclidean distance given by, for $0\leq z\leq 2$ (see Lemma 1.4.4 in
 \cite{Ivanov89})
\begin{equation}
\frac{d}{dz}P\left[ \Vert W_{1}-W_{2}\Vert \leq z\right] =\frac{1}{\sqrt{\pi
}}\Gamma \left( \frac{d}{2}\right) \Gamma ^{-1}\left( \frac{d-1}{2}\right)
z^{d-2}\left( 1-\frac{z^{2}}{4}\right) ^{\frac{d-3}{2}}.
\label{distdsphe}
\end{equation}

\medskip

\noindent \emph{Condition L}.
\begin{itemize}
\item[(i)] Assume that $\sup_{z\in [0,2]}C(z,\tau)\to 0$ as $\tau
\to \infty.$

\item[(ii)] There exists $\delta \in (0,1),$ such that
\begin{eqnarray}
&&\lim_{T\rightarrow \infty }\frac{1}{T^{\delta }}\int_{0}^{T}\left( 1-\frac{%
\tau }{T}\right) \int_{0}^{2}z^{d-2}C(z,\tau )\left( 1-\frac{z^{2}}{4}%
\right) ^{\frac{d-3}{2}}dzd\tau =\infty .  \notag  \label{csphrf}
\end{eqnarray}%

\end{itemize}
\begin{remark}
Note that, for the restriction to the sphere of a  STGRF  with covariance function as in \textbf{Example 1}, Condition L
holds if $d\geq 3.$
\end{remark}

Under \emph{Conditions K and  L}, applying the corresponding reduction
theorem,  the following central limit result follows.

\begin{theorem}
\label{thclrsrf} Assume that \emph{Conditions K} and  \emph{L} hold. Then, the
random variable
\begin{equation*}
\frac{N_{T}-T(1-\Phi (u))|\mathbb{S}_{d-1}(1)|}{\phi (u)\left[ \frac{8\pi ^{d-%
\frac{1}{2}}}{\Gamma \left( \frac{d}{2}\right) \Gamma \left( \frac{d-1}{2}%
\right) }T\int_{0}^{T}\left( 1-\frac{\tau }{T}\right) \int_{0}^{2}C(z,\tau
)z^{d-2}\left( 1-\frac{z^{2}}{4}\right) ^{\frac{d-3}{2}}dzd\tau \right]
^{1/2}}
\end{equation*}%
\noindent has asymptotically standard normal distribution as $T\rightarrow
\infty .$
\end{theorem}

The following central limit result is obtained for $T$--varying thresholds.

\begin{theorem}
\label{th2sph} Assume that \emph{Conditions} K and L hold, and there
exists  $\beta \in (0,\delta )$ such that as $T\rightarrow \infty ,$
\begin{equation*}
u^{2}(T)\sup_{u\in (0,2)}C(u,T^{\beta })\rightarrow 0,\quad u^{2}(T)=o(\log
(T)).
\end{equation*}%
Then,  the random variable
\begin{equation*}
\frac{N_{T}^{\star }-T(1-\Phi (u(T))|\mathbb{S}_{d-1}(1)|}{\phi (u(T))\left[
\frac{8\pi ^{d-\frac{1}{2}}}{\Gamma \left( \frac{d}{2}\right) \Gamma \left(
\frac{d-1}{2}\right) }T\int_{0}^{T}\left( 1-\frac{\tau }{T}\right)
\int_{0}^{2}C(z,\tau )z^{d-2}\left( 1-\frac{z^{2}}{4}\right) ^{\frac{d-3}{2}%
}dzd\tau \right] ^{1/2}},
\end{equation*}%
\noindent has asymptotically a standard normal distribution, as $T\rightarrow \infty ,$ where
\begin{equation*}
N_{T}^{\star }=\left\vert \left\{ 0\leq t\leq T;\ T_{R}(x,t)\geq u(T),\ x\in
\mathbb{S}_{d-1}(1)\right\} \right\vert .
\end{equation*}%
\end{theorem}

The proofs of  Theorems \ref{thclrsrf} and \ref{th2sph}  can be  obtained from equations (\ref{eqchaoscomp}) and (\ref{distdsphe}),
 in a similar way to the proofs of  Theorems \ref{MFth1} and \ref{th3b}, respectively.

\subsection{Spherical spatiotemporal covariance functions}
Special cases of stationary covariance functions on spheres cross
time have recently been analyzed in  \cite{White19}. In our
case, we pay special attention to the family of nonseparable covariance
functions introduced in equation (11) in \cite{White19}, since  its restriction to $\mathbb{S}_{d-1}(1)$   can be considered as proposed here. In addition,
in Theorem 2 in \cite{White19},   competitive models of spatiotemporal
spherical  covariance functions are proposed for real  data  analysis.  In particular,  the covariance function family

$$C(\theta ,u)=\frac{\sigma^{2}}{\psi(u^{2})}\varphi \left(\frac{\theta }{\psi(u^{2})}\right),\quad \theta \in [0,\pi],\ u\in \mathbb{R},$$
\noindent is considered for  surface air temperature reanalysis data.  These covariance modes  capture  the strong spatial structure  displayed by data given by daily temperature averages over a global grid.  Since the overall temperature distribution is similar across days displaying  a clear spatial structure, the implemented spatiotemporal  spherical covariance models that rescale space with time, and are expressed in terms of the geodesic  spherical distance instead of the Euclidean distance, allow to  effectively capture the strong spatial structure in this type of data.

 Note also that, one can consider the restriction to the sphere of  the covariance function family considered in equations (6) and (11) in \cite{White19}, and beyond.  Equations  (\ref{eqexorth}) and (\ref{bimeasureR}) then hold for such a  restriction. Specifically,  in  equation (\ref{gneit}), we consider, for every $u\geq 0,$

\begin{eqnarray}
\varphi (u)&=&\left( 2^{\nu -1}\Gamma \left( \nu \right) \right) ^{-1}\left(
cu^{1/2}\right) ^{\nu }K_{\nu }\left( cu^{1/2}\right) ,\text{ }c>0,\text{ }%
\nu >0\nonumber\\
\psi (u)&=&(1+au^{\alpha })^{\beta },\text{ }a>0,\text{ }0<\alpha \leq 1,\text{
}0<\beta \leq 1,\label{fcfsrdb}
\end{eqnarray}
\noindent where $K_{\nu }(z)$ is  the modified Bessel function of the
second kind of order $\nu,$ or MacDonald function (see, eg., \cite{Gradshteyn07}). Thus, \begin{equation}
\varphi(\|z\|^{2})=\sigma ^{2}\left( 2^{\nu -1}\Gamma \left( \nu \right) \right)
^{-1}\left( c\left\Vert z\right\Vert \right) ^{\nu }K_{\nu }\left(
c\left\Vert z\right\Vert \right) ,\text{ }\sigma ^{2}>0,\text{ }c>0,\text{ }%
\nu >0, \ z\in \mathbb{R}^{d} \label{wmcov}
\end{equation}
\noindent with associated  Fourier transform
\begin{eqnarray}
\widehat{\varphi }(\lambda )&=& \mathcal{M}
 \left( c^{2}+\left\Vert
\lambda \right\Vert ^{2}\right) ^{-\left( \nu +\frac{d}{2}\right) },\quad \lambda
\in \mathbb{R}^{d},\quad \mathcal{M}>0,\label{fcfsrdb3}
\end{eqnarray}

\noindent that is involved in the definition of the spectral measure $F$
in equation (\ref{sm}).  Its restriction to $\mathbb{S}_{d-1}(1)$ then leads to an alternative  family of spherical covariance functions to the ones considered in \cite{White19}.

\subsection*{Acknowledgements}
Nikolai Leonenko (NL) would like to thank for support and hospitality during the programme Fractional Differential Equations and the programme Uncertainly Quantification and Modelling of Materials in  Isaac Newton Institute for Mathematical Sciences, Cambridge

Nikolai Leonenko would like
to thank for support and hospitality during the programmes "Fractional Differential Equations" (2022), "Uncertainly Quantification and Modelling of Materials"(2024) and "Stochastic systems for anomalous diffusion"(2024) in Isaac Newton Institute for Mathematical Sciences, Cambridge.

\medskip

\noindent The first author  was partially supported by the Croatian Science Foundation (HRZZ) grant Scaling in Stochastic Models (IP-2022-10-8081).
He also
was supported in part  by ARC Discovery
Grant DP220101680 (Australia), LMS grant 42997 (UK) and grant FAPESP
22/09201-8 (Brazil).

\noindent The second author was supported in part by MCIN/ AEI/PID2022--142900NB-I00,  and CEX2020-001105-M MCIN/ AEI/10.13039/501100011033).

\end{document}